\documentclass[12pt,a4paper,equation]{article}
\usepackage[greek,american]{babel}
\usepackage[iso-8859-7]{inputenc}
\usepackage{amssymb,amsfonts}
\usepackage[dvips]{graphicx}
\usepackage{amsmath}
\usepackage{epsfig}
\usepackage{latexsym}
\usepackage{makeidx}
\usepackage{pst-math,pst-xkey}

\numberwithin{equation}{section}
\makeindex
\begin{document}
\date{}
\author{{\bf Vassilis G. Papanicolaou}\\
Department of Mathematics\\
National Technical University of Athens\\
Zografou Campus\\
157 80 Athens, GREECE\\
\underline{papanico@math.ntua.gr}}
\title{Some Results on Ordinary Differential Operators with Periodic Coefficients}
\maketitle
\begin{abstract}
For a general ordinary differential operator $\mathcal{L}$ with periodic coefficients we prove that the characteristic polynomial of the
Floquet matrix is irreducible over the field of meromorphic functions. We also consider a multipoint eigenvalue problem and show that its
eigenspaces are spanned by pure or generalized Floquet solutions. Finally, at the end of the paper we mention some relevant conjectures and open
questions.
\end{abstract}


\smallskip

\noindent\textbf{2010 AMS Mathematics Classification.} 47E05;  34B10;  47A75.

\section{Preliminaries}

The object of this article is the study of the equation
\begin{equation}
\mathcal{L} u = \lambda u,
\label{A3}
\end{equation}
where $\lambda \in \mathbb{C}$ is the spectral parameter and $\mathcal{L}$ is the differential operator
\begin{equation}
\mathcal{L} u := \frac{d^n u}{dx^n} + \sum_{k=0}^{n-1} p_k(x) \frac{d^k u}{dx^k},
\qquad\qquad
n \geq 2,
\label{A1}
\end{equation}
whose  coefficients $p_k(x)$, $k = 0, 1, \dots, n-1$, are assumed
complex-valued, infinitely differentiable, and $b\,$-periodic, namely
\begin{equation}
p_k(x+b) = p_k(x)
\qquad \text{for all }\; x \in \mathbb{R}, \quad k = 0, 1, \dots, n-1,
\label{A2}
\end{equation}
where the period $b > 0$ is fixed. We view $\mathcal{L}$ as a densely defined unbounded operator acting on $L_2(\mathbb{R})$.
Let us notice that the smoothness assumption on the coefficients of $\mathcal{L}$ is a convenience assumption and can be relaxed.

It is well known \cite{N} that via the transformation
\begin{equation*}
u(x) = v(x) \exp\left(-\frac{1}{n} \int_0^x p_{n-1}(\xi) d\xi\right)
\end{equation*}
one sees that without loss of generality we can take
\begin{equation}
p_{n-1}(x) \equiv 0.
\label{A0}
\end{equation}
From now on we will always assume (\ref{A0}).

For $n = 2$ the operator $\mathcal{L}$ becomes the so-called Hill operator
(also known as the one-dimensional Schr\"{o}dinger operator with a periodic potential). Usually, in this case $p_0(x)$ is assumed real-valued.
There is, however, considerable literature even for the case of a nonreal $p_0(x)$ (see, e.g., \cite{Tk}).

The (formal) adjoint operator $\mathcal{L}^{\ast}$ of $\mathcal{L}$ is given by
\begin{equation}
\mathcal{L}^{\ast} u := (-1)^n \frac{d^n u}{dx^n} + \sum_{k=0}^{n-2} (-1)^k \frac{d^k}{dx^k} \left[\overline{p_k(x)} \, u \right],
\label{A3a}
\end{equation}
where $\overline{p_k(x)}$ denotes the complex conjugate of $p_k(x)$. Notice that $\mathcal{L}^{\ast}$ is of the same form as $\mathcal{L}$
(except for the minus sign in front of $d^n/dx^n$ in the case where $n$ is odd).

If $\mathcal{L}$ is symmetric \cite{R-S} and has real coefficients, then it must have the form
\begin{equation}
\mathcal{A} u := \frac{d^{2\nu} u}{dx^{2\nu}}
+ \sum_{k=0}^{\nu-1} \frac{d^k}{dx^k} \left[ a_k(x) \frac{d^k u}{dx^k}\right],
\label{A1aa}
\end{equation}
where $\nu \geq 1$ ($n = 2\nu$) and the coefficients $a_k(x)$, $k = 0, 1, \dots, \nu-1$, are real-valued,
$C^{\infty}$, and $b\,$-periodic.

It is well known \cite{D-S} that $\mathcal{A}$, viewed as a symmetric operator acting on $L_2(\mathbb{R})$, has a unique self-adjoint
extension which we also denote by $\mathcal{A}$.

The ``unperturbed'' operator  associated to $\mathcal{L}$ is
\begin{equation}
\tilde{\mathcal{L}} := \frac{d^n}{dx^n},
\qquad \qquad
(\text{i.e. the $n$-th derivative}).
\label{A1a}
\end{equation}

Apart from its own significance, the operator $\tilde{\mathcal{L}}$ plays an important role in the analysis of $\mathcal{L}$ since the
asymptotic behavior, as $\lambda \to \infty$, of many quantities related to
(\ref{A3}) is determined by the associated quantities corresponding to the
unperturbed case $\tilde{\mathcal{L}} u = \lambda u$ (usually, the latter quantities can be computed explicitly).
Furthermore, $\mathcal{L}$ and (the solutions of) (\ref{A3}) can be considered as perturbations or continuous deformations of
$\tilde{\mathcal{L}}$ and $\tilde{\mathcal{L}} u = \lambda u$ respectively. From now on, we will be using the notational convention that tilded
quantities correspond to the unperturbed case.

\subsection{Some basic facts}

Consider an one-parameter family of operators of the form (\ref{A1}) with $p_{n-1} \equiv 0$, namely
\begin{equation}
\mathcal{L}(t) u := \frac{d^n u}{dx^n} + \sum_{k=0}^{n-2} p_k(x; t) \frac{d^k u}{dx^k},
\qquad
t \in \mathbb{T},
\label{TA1}
\end{equation}
where, typically, $\mathbb{T}$ is an interval of $\mathbb{R}$. For example, we can have $\mathbb{T} = [0, 1]$ with
$\mathcal{L}(0) = \tilde{\mathcal{L}}$ and $\mathcal{L}(1) = \mathcal{L}$, so that, intuitively, $\mathcal{L}$ can be seen as a deformation
of $\tilde{\mathcal{L}}$. We always assume that each $p_k(x; t)$ depends smoothly on $t$ and that there is a constant $M > 0$ such that
\begin{equation}
|p_k(x; t)| \leq M
\qquad\text{for all }\;
x \in \mathbb{R}, \ \ t \in \mathbb{T}, \ \ k = 0, 1, \dots, n-2.
\label{TA111}
\end{equation}
A special such family of operators appears in Subsection 2.5.

Let us also consider the associated differential equation
\begin{equation}
\mathcal{L}(t) u = \lambda u,
\qquad
\lambda \in \mathbb{C}.
\label{TA3}
\end{equation}
The proof of the following important bounds regarding certain solutions of (\ref{TA3}) can be found in \cite{N} (Part I, Ch. II).

Put
\begin{equation}
\lambda = \zeta^n
\label{TA2}
\end{equation}
and divide the complex $\zeta$-plane into the $2n$ sectors
\begin{equation}
S_l := \left\{\zeta \in \mathbb{C}:\, \frac{l \pi}{n} \leq \arg(\zeta) \leq \frac{(l+1) \pi}{n} \right\},
\quad
l = 0, 1, \dots, 2n-1
\label{TA4}
\end{equation}
(notice that adjacent sectors have a ray in common and also that we can cover all the values of $\lambda \in \mathbb{C}$ by taking
$\zeta$ to be the principal $n$-th root of $\lambda$, i.e. $\zeta \in S_0 \cup S_1^{\circ}$, where $S_1^{\circ}$ denotes the
interior of $S_1$). Furthermore, let
\begin{equation}
0 \leq x \leq B < \infty,
\label{TA5}
\end{equation}
where $B > 0$ is a fixed number.

Then, for each region $S_l$ of the complex $\zeta$-plane and for $|\zeta|$ sufficiently large, say $|\zeta| > Z$, where $Z$ is a
positive constant, which can be taken independent of $t$,
the equation (\ref{TA3}) possesses $n$ linearly independent solutions
$y_j = y_j(x; \zeta; t)$, $j = 1, \dots, n$ (these solutions are not unique), satisfying
(for all $x \in [0, B]$ and $t \in \mathbb{T}$)
\begin{equation}
\left|\frac{y_j^{(k-1)}(x; \zeta; t)}{\omega_j^{k-1} \zeta^{k-1} e^{\omega_j \zeta x}} - 1 \right| \leq \frac{K}{|\zeta|},
\qquad
j, k = 1, \dots, n,
\label{TA6}
\end{equation}
where $\omega_1, \dots, \omega_n$ are the $n$-th roots of unity, $K$ is a constant which is independent of $t$, $x$, and $\zeta$
(i.e. $K$ depends only on the coefficients $p_0(x), \dots, p_{n-2}(x)$ of $\mathcal{L} = \mathcal{L}(1)$ and the constant $B$), and
$y_j^{(k-1)}(x; \zeta; t)$ denotes the $(k-1)$-derivative of $y_j$ with respect to $x$. Notice that each $\omega_j$, $j = 1, \dots, n$,
is associated to a different $y_j$ and, also, that the $y_j$'s are linearly independent as functions of $x$, for each fixed pair $(\zeta, t)$
(recall that $|\zeta| > Z$). Furthermore, for each $l = 0, 1, \dots, 2n-1$, every function $y_j^{(k-1)}(x; \zeta; t)$, $j,k = 1, \dots, n$
is analytic in $\zeta \in S_l$, $|\zeta| > Z$.
We refer to such solutions $y_j$ as \textit{Naimark solutions} of (\ref{TA3}) and, in particular (when $t=1$), of (\ref{A3}).
Of course, (\ref{TA6}) is essentially an asymptotic formula since it is equivalent to
\begin{equation}
y_j^{(k-1)}(x; \zeta; t) = \omega_j^{k-1} \zeta^{k-1} e^{\omega_j \zeta x} \left[1 + O\left(\frac{1}{\zeta}\right) \right],
\qquad
j, k = 1, \dots, n,
\label{TA6a}
\end{equation}
as $\zeta \to \infty$.

In the unperturbed case ($t = 0$) we can take
\begin{equation}
\tilde{y}_j(x; \zeta) := e^{\omega_j \zeta x},
\qquad
j = 1, \dots, n.
\label{TA7}
\end{equation}

For $\zeta \in S_l$ (for any $l$) with $|\zeta| > Z$, the Wronskian $W_y := W[y_1, \dots, y_n]$ of $y_1, \dots, y_n$ is independent of $x$
due to (\ref{A0}), hence it is determined by its value at $x = 0$. Using (\ref{TA6}) with $x = 0$ yields
\begin{equation}
W_y = \zeta^{n(n-1)/2}
\left|
\begin{array}{cccc}
1 + O(\zeta^{-1}) & 1 + O(\zeta^{-1}) & \cdots & 1 + O(\zeta^{-1})  \\
\omega_1[1 + O(\zeta^{-1})] & \omega_2[1 + O(\zeta^{-1})] & \cdots & \omega_n[1 + O(\zeta^{-1})] \\
\vdots & \vdots & \ddots & \vdots  \\
\omega_1^{n-1}[1 + O(\zeta^{-1})] & \omega_2^{n-1}[1 + O(\zeta^{-1})] & \cdots & \omega_n^{n-1}[1 + O(\zeta^{-1})]
\end{array}
\right|
\label{TA7a}
\end{equation}
as $\zeta \to \infty$, $\zeta \in S_l$. Now, since
\begin{equation}
\{\omega_1, \omega_2, \dots, \omega_n\} = \{1, \rho, \dots, \rho^{n-1}\},
\qquad \text{where }\;
\rho := e^{2\pi i / n},
\label{A6b}
\end{equation}
the determinant in the right-hand side of (\ref{TA7a}) is equal to
\begin{equation}
\pm \left|
\begin{array}{cccc}
1 + O(\zeta^{-1}) & 1 + O(\zeta^{-1}) & \cdots & 1 + O(\zeta^{-1})  \\
1 + O(\zeta^{-1}) & \rho [1 + O(\zeta^{-1})] & \cdots & \rho^{n-1} [1 + O(\zeta^{-1})] \\
\vdots & \vdots & \ddots & \vdots  \\
1 + O(\zeta^{-1}) & \rho^{n-1}[1 + O(\zeta^{-1})] & \cdots & \rho^{(n-1)^2}[1 + O(\zeta^{-1})]
\end{array}
\right| .
\label{TA7b}
\end{equation}
By expanding the determinant appearing in (\ref{TA7b}) using the elements of the $n$-th row and their associated $(n-1) \times (n-1)$ minors,
and then applying straightforward induction on $n$ we can deduce that the determinant is equal to $R_n [1 + O(\zeta^{-1})]$, where
\begin{equation}
R_n := \left|
\begin{array}{cccc}
1 & 1 & \cdots & 1  \\
1 & \rho & \cdots & \rho^{n-1} \\
\vdots & \vdots & \ddots & \vdots  \\
1 & \rho^{n-1} & \cdots & \rho^{(n-1)^2}
\end{array}
\right|
= \varepsilon_n \, n^{n/2},
\quad \text{with }\;
\varepsilon_n \in \{1, -1, i, -i\}
\label{TA7c}
\end{equation}
(the second equality follows, e.g., by observing that $R_n^2$ is the discriminant of the polynomial $z^n - 1$,
hence $R_n^2 = (-1)^{(n+2)(n-1)/2} n^n$). Therefore, (\ref{TA7a}) becomes
\begin{equation}
W_y = W[y_1, \dots, y_n] = \pm R_n \lambda^{(n-1)/2} \left[1 + O\left(\frac{1}{\lambda^{1/n}}\right)\right],
\qquad
\lambda \to \infty,
\label{TA7d}
\end{equation}
uniformly in $t$. Notice that for the unperturbed case we have $W_{\tilde{y}} = \pm R_n \lambda^{(n-1)/2}$.

\subsection{The fundamental solutions}

Let ${\cal V}_{\lambda}$ be the $n$-dimensional vector space consisting of the solutions of (\ref{A3}) (or of (\ref{TA3}) for a fixed $t$).
We introduce the standard basis of ${\cal V}_{\lambda}$, namely the fundamental solutions $u_j(x; \lambda)$:
\begin{equation}
u_j^{(i-1)}(0; \lambda) = \delta_{ij},
\qquad
i,j = 1, \dots, n,
\label{A5}
\end{equation}
where $u_j^{(i-1)}$ denotes the $(i-1)$-th derivative of $u_j$ with respect to $x$ and $\delta_{ij}$ is the Kronecker delta.
For a fixed $x > 0$ each $u_j(x; \lambda)$, $j = 1, \dots, n$, as well as its first $n$ derivatives with respect to $x$,
are entire functions of $\lambda$, whose order can be shown by standard arguments to be $1 / n$ (and the associated type is $x$).

Using this standard basis, the solution
\begin{equation*}
u(x) = c_1 u_1(x; \lambda) + c_2 u_2(x; \lambda) + \cdots + c_n u_n(x; \lambda)
\end{equation*}
of (\ref{A3}) (or (\ref{TA3}), i.e. the (unique) solution $u(x)$ such that $u^{(j-1)}(0) = c_j$ for $j = 1, \dots, n$,
is represented by the vector
\begin{equation*}
u = [c_1, c_2, \dots, c_n]^{\top}.
\end{equation*}

In the unperturbed case the fundamental solutions are
\begin{align}
\tilde{u}_1(x; \lambda) &= \frac{1}{n} \left[
\exp(\lambda^{1 / n} x) +  \exp(\rho \lambda^{1 / n} x) + \cdots + \exp\left(\rho^{n-1} \lambda^{1 / n} x\right)\right]
\nonumber
\\
&= \sum_{k=0}^{\infty}\, \frac{x^{nk}}{(nk)!} \, \lambda^k
\label{A6a}
\end{align}
($\lambda^{1 / n}$ denotes the principal branch of the $n$-th root function, i.e. $0 \leq \arg(\lambda^{1 / n}) < 2\pi / n$ and $\rho$
is given by (\ref{A6b})) and
\begin{equation}
\tilde{u}_j(x; \lambda) = \frac{\tilde{u}_1^{(n-j+1)}(x; \lambda)}{\lambda},
\qquad \qquad \qquad \qquad
j = 2, \dots, n. \qquad
\label{A6aa}
\end{equation}

\smallskip

\textbf{Remark 1.} A common special case of an old result, that goes back to Laguerre, states the following: Suppose $h(z)$ is an
entire function of genus $0$ or $1$, which is real for real $z$, and all its zeros are real and negative (for the definition of the
genus of an entire function see, e.g., \cite{A}, \cite{H}, or \cite{T}). Then
\begin{equation}
f(\lambda) := \sum_{k=0}^{\infty}\, \frac{h(k)}{k!} \, \lambda^k
\label{A6d}
\end{equation}
is entire and all its zeros are real and negative. The proof of this, as well as the more general result,
can be found in \cite{H} or \cite{T}. If we take
\begin{equation}
h(z) := \frac{\Gamma(z+1)}{\Gamma(nz+1)}\, e^{n z \ln x},
\label{A6e}
\end{equation}
where $\Gamma(\cdot)$ denotes the Gamma function, then $f(\lambda)$ of (\ref{A6d}) becomes $\tilde{u}_1(x; \lambda)$ of (\ref{A6a}). Thus
all zeros of $\tilde{u}_1(x; \lambda)$, viewed as a function of $\lambda$, are
real and negative for any fixed $x > 0$. Furthermore, since the order of $\tilde{u}_1(x; \lambda)$ is $1/n$, the Hadamard Factorization Theorem
implies that there are infinitely many such zeros. In the same way one can show that all zeros of $\tilde{u}_j(x; \lambda)$ (as a function of $\lambda$) are real and negative, for any fixed $x > 0$ and for any $j = 1, \dots, n$.

On the other hand, if we fix a $\lambda < 0$,
then, by a ``shooting" argument (i.e. starting with an $x$ very close to 0, so that $\tilde{u}_j(x) > 0$, and then moving $x$ continuously
away from $0$, since, initially, $\tilde{u}^{(n)}_j(x) = \lambda \tilde{u}_j(x) < 0$, at some value, say, $x = x^j_1$ we will have
$\tilde{u}_j(x^j_1) = 0$; as $x$ gets bigger than $x^j_1$ we will have $\tilde{u}_j(x) < 0$ until $x$ reaches a value, say $x = x^j_2$, where
$\tilde{u}_j(x^j_2) = 0$, and so on)
it is not hard to see that each $\tilde{u}_j(x; \lambda)$, $j = 1, \dots, n$ (viewed now as a function of $x$), has infinitely many zeros
$x^j_1 < x^j_2 < \cdots$ in $(0, \infty)$, which cannot accumulate at a finite point since $\tilde{u}_j(x; \lambda)$ is also entire in $x$,
for every $j = 1, \dots, n$. These zeros interlace in the sense
\begin{equation}
x^n_1 < x^{n-1}_1 < \cdots < x^1_1 < x^n_2 < x^{n-1}_2 < \cdots < x^1_2 < x^n_3 < \cdots
\label{A6f}
\end{equation}
and, in view of (\ref{A6aa}), they are simple zeros.
\hfill $\diamondsuit$

\smallskip

(The symbol $\diamondsuit$ indicates the end of a remark).

\smallskip

Let us now return to the fundamental solutions of (\ref{A3}) (or of (\ref{TA3}) for a fixed $t$). Assuming $\zeta \in S_l$, for some $l$,
and $|\zeta| > Z$ the $1$st fundamental solution can be written as
\begin{equation}
u_1(x; \lambda) = \frac{1}{W_y} \left|
\begin{array}{cccc}
y_1(x) & y_2(x) & \cdots & y_n(x)  \\
y_1'(0) & y_2'(0) & \cdots & y_n'(0) \\
\vdots & \vdots & \ddots & \vdots  \\
y_1^{(n-1)}(0) & y_2^{(n-1)}(0) & \cdots & y_n^{(n-1)}(0)
\end{array}
\right|,
\label{A6g}
\end{equation}
where $y_1, \dots, y_n$ are a choice of Naimark solutions satisfying (\ref{TA6}) and $W_y$ is their Wronskian. Likewise,
\begin{equation}
u_j(x; \lambda) = \frac{(-1)^{j-1}}{W_y} \det \left[
\begin{array}{c}
\vec{y}(x)  \\
M_j
\end{array}
\right],
\label{A6h}
\end{equation}
where $\vec{y}(x)$ is the row vector
\begin{equation*}
\vec{y}(x) := \left[y_1(x), \dots, y_n(x) \right]
%\label{A6i}
\end{equation*}
and $M_j$ is the $(n-1) \times n$ matrix obtained by the matrix
\begin{equation*}
M := \left[
\begin{array}{cccc}
y_1(0) & y_2(0) & \cdots & y_n(0)  \\
y_1'(0) & y_2'(0) & \cdots & y_n'(0) \\
\vdots & \vdots & \ddots & \vdots  \\
y_1^{(n-1)}(0) & y_2^{(n-1)}(0) & \cdots & y_n^{(n-1)}(0)
\end{array}
\right],
%\label{A6j}
\end{equation*}
after erasing its $j$-th row.

Using (\ref{TA6}) in (\ref{A6g}) yields (as $\zeta \to \infty$, $\zeta \in S_l$)
\begin{equation}
u_1(x; \lambda) = \frac{\zeta^{n(n-1)/2}}{W_y} \left|
\begin{array}{cccc}
e^{\omega_1 \zeta x} [1 + O(\zeta^{-1})] & e^{\omega_2 \zeta x} [1 + O(\zeta^{-1})] & \cdots & e^{\omega_1 \zeta x} [1 + O(\zeta^{-1})]  \\
\omega_1 [1 + O(\zeta^{-1})] & \omega_2 [1 + O(\zeta^{-1})] & \cdots & \omega_n [1 + O(\zeta^{-1})] \\
\vdots & \vdots & \ddots & \vdots  \\
\omega_1^{n-1} [1 + O(\zeta^{-1})] & \omega_2^{n-1} [1 + O(\zeta^{-1})] & \cdots & \omega_n^{n-1} [1 + O(\zeta^{-1})]
\end{array}
\right|
\label{A6k}
\end{equation}
(in the unperturbed case formula (\ref{A6k}) is exact, i.e. the terms $O(\zeta^{-1})$ are all identically $0$).
In view of (\ref{A6a}), (\ref{TA7d}), and Remark 1, formula (\ref{A6k}) implies (by expanding the determinant of (\ref{A6k}) using the
elements of the first row and their associated $(n-1) \times (n-1)$ minors) that if $\lambda$ approaches $\infty$ along a ray
other than the negative real axis (i.e. $\arg(\lambda)$ is fixed and $\ne \pi$), then
\begin{equation}
u_1(x; \lambda) = \tilde{u}_1(x; \lambda) \left[1 + O\left(\frac{1}{\lambda^{1/n}}\right)\right],
\label{A6l}
\end{equation}
uniformly in $x \in [0, B]$ and $t \in \mathbb{T}$. In the same way one can show that (as $\lambda \to \infty$ along a fixed ray
other than the negative real axis)
\begin{equation}
u_j(x; \lambda) = \tilde{u}_j(x; \lambda) \left[1 + O\left(\frac{1}{\lambda^{1/n}}\right)\right]
\qquad \text{for all }\;
j = 1, \dots, n,
\label{A6m}
\end{equation}
uniformly in $x \in [0, B]$ and $t \in \mathbb{T}$.

\section{Floquet theory}

Let $\mathcal{T}$ be the
$b\,$-shift (linear) operator defined by
\begin{equation}
(\mathcal{T}f)(x) := f(x+b).
\label{A4}
\end{equation}
Then, it is clear that $\mathcal{T}$ maps solutions of (\ref{A3}) to solutions of
(\ref{A3}), i.e. $\mathcal{T}$ and $\mathcal{L}$ commute. It follows that, for each $\lambda \in \mathbb{C}$,
the $b\,$-shift $\mathcal{T}$ can be viewed as an operator acting on the
$n$-dimensional vector space ${\cal V}_{\lambda}$ of the solutions of (\ref{A3}). This simple
observation is the essence of Floquet theory.

Let us consider the standard basis of ${\cal V}_{\lambda}$, namely the fundamental solutions $u_j(x; \lambda)$ of (\ref{A3})
(see (\ref{A5})). The matrix of $\mathcal{T}$ with respect to this basis is
\begin{equation}
T = T(\lambda)
= \left[t_{ij}(\lambda)\right]_{1 \leq i,j \leq n}
\qquad \text{with }\;
t_{ij}(\lambda) = u_j^{(i-1)}(b; \lambda), \quad 1 \leq i,j \leq n.
\label{A6}
\end{equation}
This is the \textit{Floquet matrix} of (\ref{A3}) and, as we have seen,
its elements $u_j^{(i-1)}(b; \lambda)$ are entire functions of $\lambda$, whose order is $1 / n$ (and the associated type is $b$).
Notice, also, that (\ref{A0}) implies
\begin{equation}
\det T(\lambda) = 1.
\label{A6c}
\end{equation}

\textbf{Theorem 1.} Let $\lambda = \zeta^n$, where $\zeta \in S_l$ for some $l = 0, 1, \dots, 2n-1$ (see (\ref{TA4})) and $|\zeta| > Z$.
If $Y := [\,\mathcal{T}\,]_y$, namely if $Y$ is the matrix of the $b\,$-shift operator $\mathcal{T}$ with respect to the basis
$y_1, \dots, y_n$ of a choice of Naimark solutions of (\ref{A3}) (or (\ref{TA3})) satisfying (\ref{TA6}), then
\begin{equation}
Y = \left[\gamma_{ij} \right]_{1 \leq i,j \leq n},
\quad \text{where }\;
\gamma_{ij} = e^{\omega_j \zeta b}\left[\delta_{ij} + O(\zeta^{-1}) \right],
\quad
\zeta \to \infty
\ \; (\zeta \in S_l)
\label{TA8}
\end{equation}
(again, $\delta_{ij}$ denotes the Kronecker delta).

\smallskip

\textit{Proof}. Fix an $l$ in $\{0, \dots, (2n-1)\}$ and assume $\zeta \in S_l$ and $|\zeta| > Z$. Since
$(\mathcal{T} y_j)(x) = y_j(x+b)$ the elements
$\gamma_{1j}, \dots, \gamma_{nj}$ of the $j$-th column of $Y$ must satisfy (for each $j = 1, \dots, n$)
\begin{equation}
y_j(x+b) = \gamma_{1j} \, y_1(x) + \cdots + \gamma_{nj} \, y_n(x).
\label{TA9}
\end{equation}
Differentiating (\ref{TA9}) with respect to $x$ repeatedly yields
\begin{equation}
y_j^{(k)}(x+b) = \gamma_{1j} \, y_1^{(k)}(x) + \cdots + \gamma_{nj} \, y_n^{(k)}(x),
\qquad
k = 0, \dots, n-1.
\label{TA10}
\end{equation}
If we set $x = 0$ in (\ref{TA10}) and invoke (\ref{TA6}) (where $B$ can be any fixed number greater than $b$) we get, as $\zeta \to \infty$,
\begin{equation}
\omega_j^k e^{\omega_j \zeta b} \left[1 + O(\zeta^{-1}) \right]
= \omega_1^k \left[1 + O(\zeta^{-1}) \right] \gamma_{1j} + \cdots \, + \omega_n^k \left[1 + O(\zeta^{-1}) \right] \gamma_{nj},
\label{TA11}
\end{equation}
or, by substituting $\gamma_{ij} = e^{\omega_j \zeta b} \, \hat{\gamma}_{ij}$, $i = 1, \dots, n$, and elimining one of the factors of
the form $\left[1 + O(\zeta^{-1})\right]$,
\begin{equation}
\omega_1^k \left[1 + O(\zeta^{-1}) \right] \hat{\gamma}_{1j} + \cdots \, + \omega_n^k \left[1 + O(\zeta^{-1}) \right] \hat{\gamma}_{nj}
= \omega_j^k,
\quad
k = 0, \dots, n-1,
\label{TA11a}
\end{equation}
as $\zeta \to \infty$. We can now solve the system (\ref{TA11a}) for $\hat{\gamma}_{ij}$, $i = 1, \dots, n$, by, say, Cramer's rule.
Noticing that, in (\ref{TA11a}) the determinant of the coefficients of the $\hat{\gamma}_{ij}$, $i = 1, \dots, n$, is equal to
$\pm R_n \left[1 + O(\zeta^{-1}) \right]$, where $R_n$ is given by (\ref{TA7c}), we obtain (\ref{TA8}).
\hfill $\blacksquare$

\smallskip

(As usual, the symbol $\blacksquare$ indicates the end of a proof).

\smallskip

Observe that, in the $t$-dependent case (\ref{TA3}) the asymptotic formula (\ref{TA8}) is uniform in $t \in \mathbb{T}$, since the constant
$K$ of (\ref{TA6}) is independent of $t$.

In view of (\ref{TA7}), formula (\ref{TA8}) implies that there is a matrix $E = \left[\epsilon_{ij} \right]_{1 \leq i,j \leq n}$ such that
\begin{equation}
Y = (I + E) \tilde{Y},
\qquad \text{with }\;
E = O(\zeta^{-1})
\quad \text{as }\;
\zeta \to \infty
\quad (\zeta \in S_l)
\label{TA12}
\end{equation}
for $l = 1, \dots, 2n-1$ ($I$ is, of course, the $n \times n$ identity matrix), where $\tilde{Y}$ is the diagonal matrix
\begin{equation}
\tilde{Y} = \text{diagonal} \left(e^{\omega_1 \zeta b}, e^{\omega_2 \zeta b}, \dots, e^{\omega_n \zeta b}\right),
\label{TA12a}
\end{equation}
while by $E = O(\zeta^{-1})$ we mean that all elements $\epsilon_{ij}$ of the matrix $E$ are $O(\zeta^{-1})$.
Furthermore, $Y$, being the matrix of $\mathcal{T}$ with respect to $y_1, \dots, y_n$, is similar to $T$; in particular (see (\ref{A6c}))
$\det Y = \det T = 1$. Also, since $\omega_1, \dots, \omega_n$ are the $n$-th roots of $1$ we must have
\begin{equation}
\omega_1 + \cdots + \omega_n = 0,
\label{TA12b}
\end{equation}
which implies
\begin{equation}
\det \tilde{Y} = 1
\label{TA12c}
\end{equation}
and hence, by (\ref{TA12}) we get
\begin{equation}
\det (I + E) = 1.
\label{TA12d}
\end{equation}
%Let us close this discussion with some additional observations. By (\ref{A4}) the action of the operator $\mathcal{T}^2$ on a function $f$ is
%\begin{equation}
%(\mathcal{T}^2 f)(x) = f(x + 2b).
%\label{TA12e}
%\end{equation}
%The matrix $[\,\mathcal{T}^2\,]_y$ of $\mathcal{T}^2$ with respect to the basis $y_1, \dots, y_n$ of solutions of (\ref{A3})
%(or (\ref{TA3})) satisfying (\ref{TA6}) (for some $B > 2b$) is $Y^2$, where $Y$ is the matrix introduced in Theorem 1.
%On the other hand, if we view $2b$ as the period of the coefficients of $\mathcal{L}$, direct application of Theorem 1 implies
%that there is a matrix $E_2$ such that
%\begin{equation}
%[\,\mathcal{T}^2\,]_y = (I + E_2) \tilde{Y}^2,
%\qquad \text{with }\;
%E_2 = O(\zeta^{-1})
%\quad \text{as }\;
%\zeta \to \infty
%\quad (\zeta \in S_l),
%\label{TA12f}
%\end{equation}
%where $\tilde{Y}$ is given by (\ref{TA12a}). Thus, $Y^2 = (I + E_2) \tilde{Y}^2$ or by invoking (\ref{TA12})
%\begin{equation}
%(I + E) \tilde{Y} (I + E) = (I + E_2) \tilde{Y},
%\label{TA12g}
%\end{equation}
%which implies
%\begin{equation}
%\tilde{Y} (I + E) = (I + E') \tilde{Y},
%\qquad
%E' = O(\zeta^{-1})
%\quad \text{as }\;
%\zeta \to \infty
%\quad (\zeta \in S_l).
%\label{TA12h}
%\end{equation}

\subsection{Floquet multipliers}

The characteristic polynomial of $T$ is
\begin{equation}
P(r) = P(r; \lambda) = \det [\, r I - T(\lambda)],
\label{A7}
\end{equation}
where $I$ is the $n \times n$ identity matrix. Notice that $P(r; \lambda)$ is a monic polynomial in $r$ of degree $n$, whose coefficients are
entire functions of $\lambda$ (recall that ``monic" means that the coefficient of $r^n$ is $1$). Furthermore, (\ref{A6c}) implies that
\begin{equation}
P(0; \lambda) = \det \left[ -T(\lambda) \right] = (-1)^n.
\label{A7a}
\end{equation}
For example, in the case $n = 2$ we have $P(r; \lambda) = r^2 - [u_1 (b; \lambda) + u_2'(b; \lambda)] \, r  \, + 1$.

The eigenvalues of $T$, namely the zeros of $P(r; \lambda)$, are the \textit{Floquet multipliers} of (\ref{A3}).

Incidentally, the Floquet multipliers can be used to characterize the $L_2(\mathbb{R})$-spectrum $\sigma(\mathcal{L})$ of $\mathcal{L}$ as \cite{D-S}
\begin{equation}
\sigma(\mathcal{L}) =
\{ \lambda \, : \; |r_j(\lambda)| = 1, \ \text{for some }\; j \}.
\label{A13}
\end{equation}
Using this characterization one can show \cite{D-S} that if $\mathcal{L}$ is (essentially) self-adjoint, $\sigma(\mathcal{L})$ is a countable union
of (nondegenerate) closed intervals of the real line.
For example, it is well known that if $n = 2\nu$, then $\sigma(\tilde{\mathcal{L}})$ (recall (\ref{A1a})) is the semiaxis $[0, (-1)^\nu\infty)$,
while if $n$ is odd, then $i \tilde{\mathcal{L}}$ is essentially self-adjoint and $\sigma(i\tilde{\mathcal{L}}) = \mathbb{R}$. For a more
detailed characterization of the spectrum of self-adjoint operators with real and periodic coefficients
see, e.g., \cite{D-S}, \cite{E}, \cite{Y-S}, and \cite{K1}.

In the unperturbed case we have
\begin{equation}
\tilde{r}_j(\lambda) = \exp\left(\rho^{j-1} \lambda^{1 / n} b\right),
\qquad
j = 1, 2, \dots, n,
\label{A9b}
\end{equation}
where $\lambda^{1 / n}$ is the principal branch of the $n$-th root and $\rho$ is given by (\ref{A6b}). Notice that
$\tilde{r}_1(0) = \cdots = \tilde{r}_n(0) = 1$, i.e. for $\lambda = 0$ we have that $r = 1$ is a root of $\tilde{P}(r)$ of
multiplicity $n$ (hence, $\tilde{P}(r; 0) = (r - 1)^n$). If $\lambda \ne 0$, the multiplicity of a root of $\tilde{P}(r)$
cannot exceed $2$, since the quantities $\rho^{j-1} \lambda^{1 / n} b$, $j = 1, \dots, n$ (see (\ref{A9b})), are the vertices of
a regular $n$-gon of the complex $\lambda$-plane and, hence, no more than two of them can have the same real part
(however, let us point out that for a given $\lambda$, $\tilde{P}(r)$ may have one or several double roots). If $\zeta \ne 0$ satisfies
\begin{equation}
\zeta = \frac{2 \pi i m}{(\rho^{k-1} - \rho^{j-1}) b}
\qquad \text{for some }\;
m \in \mathbb{Z}, \ \; j,k \in \{1,\dots, n\}, \ \; j \ne k,
\label{A9c}
\end{equation}
then, for $\lambda = \zeta^n$ we must have $\tilde{r}_j(\lambda) = \tilde{r}_k(\lambda)$.
For any value of $\lambda \in \mathbb{C}$, other than $0$ and the values given by (\ref{A9c}), the polynomial $\tilde{P}(r)$
has $n$ distinct roots. Moreover, if we apply Corollary A1 of the Appendix to (\ref{A9b}), we get the following immediate consequence:

\smallskip

\textbf{Corollary 1.} For any $\beta \in (0, 1)$ there is a $\delta = \delta(\beta) > 0$, with $\delta \to 0^+$ as $\beta \to 0^+$, such that
\begin{equation}
\left| \tilde{r}_k(\lambda) - \tilde{r}_j(\lambda) \right|
\geq \beta \, \max \left\{ \left| \tilde{r}_j(\lambda) \right|, \left|\tilde{r}_k(\lambda) \right| \right\},
\label{A9e}
\end{equation}
provided
\begin{equation}
\left(\rho^{k-1} - \rho^{j-1}\right) \zeta b
\not\in \bigcup_{m \in \mathbb{Z}} \{z \in \mathbb{C}:\, |z - 2 \pi i m| < \delta\},
\label{A9f}
\end{equation}
Formula (\ref{A9f}) can be written equivalently as
\begin{equation}
\zeta \not\in E_{\delta} := \bigcup_{m \in \mathbb{Z}} \; \bigcup_{1 \leq j < k \leq n}
\left\{z \in \mathbb{C}:\, \left|z - \frac{\pi m}{b \sin\left(\frac{\pi (k-j)}{n}\right)} \, e^{-\pi(j + k - 2) i / n}\right| < \delta \right\}
\label{A9g}
\end{equation}
(by ``recycling", the $\delta$ appearing in formula (\ref{A9g}) is the $\delta$ of (\ref{A9f}) divided by $2 b \sin(\pi/n)$).

\smallskip

\textbf{Remark 2.} The set $E_{\delta}$ defined in (\ref{A9g}) is the union of countably infinitely many open disks.
Some of these disks may overlap, no matter how small $\delta$ is. However, it is clear that if $\delta$ is sufficiently small, then there is
an $\epsilon > 0$ and a strictly increasing sequence $\{a_n\}_{n \in \mathbb{N}}$ of positive numbers with $a_n \to \infty$ such that
\begin{equation}
E_{\delta} \cap \bigcup_{n \in \mathbb{N}} \left\{z \in \mathbb{C}:\, a_n \leq |z| \leq a_n + \epsilon \right\} = \emptyset.
\label{A9gh}
\end{equation}
\ \hfill $\diamondsuit$

\smallskip

Coming back to (\ref{A9c}) let us observe that, without loss of generality, we can assume that $m > 0$, i.e. we can take
$m \in \mathbb{N} := \{1, 2, \dots \}$, since we can exchange the roles of $j$ and $k$. Then, from (\ref{A9c}) it follows that the values of
$\lambda \in \mathbb{C} \setminus \{0\}$ for which $\tilde{P}(r; \lambda)$ has double roots are
\begin{equation}
\tilde{\lambda}_{m,d} := \frac{(-1)^d \pi^n m^n}{\sin(\pi d / n)^n b^n},
\qquad
m \in \mathbb{N}, \ \; d \in \{1,\dots, n-1\}
\label{A9d}
\end{equation}
(the $n$-th roots of the $\tilde{\lambda}_{m,d}$'s are the centers of the disks whose union is $E_{\delta}$).
In particular, $\tilde{\lambda}_{m,d} \in \mathbb{R}$ for all $(m, d) \in \mathbb{N} \times \{1,\dots, n-1\}$. Finally, notice that if
$n$ is odd, then $\tilde{\lambda}_{m,n-d} = -\tilde{\lambda}_{m,d}$, while if $n$ is even, then
$\tilde{\lambda}_{m,n-d} = \tilde{\lambda}_{m,d}$, hence in this case it is enough to take $d \in \{1,\dots, n/2\}$. Furthermore,
in the case of an even $n$, the value $\lambda = \tilde{\lambda}_{m,n/2}$ gives
$\tilde{r}_j(\tilde{\lambda}_{m,n/2}) = \tilde{r}_k(\tilde{\lambda}_{m,n/2}) = \pm 1$, i.e.
$\tilde{\lambda}_{m,n/2}$ is a periodic or antiperiodic eigenvalue of $\tilde{\mathcal{L}}$.

We continue with a corollary of Theorem 1, which describes the (leading) asymptotic behavior of the Floquet multipliers of (\ref{A3}) (and of (\ref{TA3})) as $\lambda \to \infty$.

\smallskip

\textbf{Corollary 2.} For each $j = 1, \dots, n$ there is an $n$-th root of unity $\omega_k$ such that the Floquet multiplier $r_j(\lambda)$
of (\ref{A3}) (and of (\ref{TA3})) satisfies
\begin{equation}
\left| \frac{r_j(\lambda)}{e^{\omega_k \zeta b}} - 1 \right|
\leq \frac{K}{|\zeta|},
\label{A10}
\end{equation}
where $\lambda = \zeta^n$ (as usual) and $K > 0$ is a constant which is independent of $\zeta$
(and of $t$ in the case of (\ref{TA3})).
Furthermore, if $|\zeta|$ is sufficiently large and $\zeta \not\in E_{\delta}$ (see (\ref{A9g})) for a sufficiently small $\delta > 0$,
then to each $r_j$ in (\ref{A10}) corresponds a different $\omega_k$.

\smallskip

\textit{Proof}. Let us fix an $l$ in $\{0, \dots, 2n-1\}$ and assume $\zeta \in S_l$, $|\zeta| > Z$.
Notice that $P(r; \lambda)$ is, also, the characteristic polynomial of the matrix $Y$ of Theorem 1, since $Y$ and $T$ are similar matrices.
We, then, write (\ref{TA12}) as $Y = \tilde{Y} + E \tilde{Y}$ and apply the Gershgorin Circle Theorem (Theorem A1 of the Appendix) to
the matrix $Y^{\top}$ in order to get (\ref{A10}).

The last statement of the corollary follows from the fact that, due to (\ref{A9e}) of Corollary 1, the Gershgorin disks
(regarding the matrix $Y^{\top}$) are disjoint for $|\zeta|$ sufficiently large, with $\zeta \not\in E_{\delta}$,
and sufficiently small $\delta > 0$.
\hfill $\blacksquare$

\smallskip

The discriminant \cite{L} $D_P(\lambda)$ of the polynomial $P(r; \lambda)$ is an
entire function of $\lambda$ and, of course, $P(r; \lambda^{\star})$ has a multiple zero if and only if
$D_P(\lambda^{\star}) = 0$. Being entire, $D_P(\lambda)$ has at most countably
many zeros, or $D_P(\lambda) \equiv 0$. However, if $D_P(\lambda) \equiv 0$, then $P(r; \lambda)$ would have multiple zeros for all
$\lambda \in \mathbb{C}$, but this is impossible since by Corrolary 1 and (\ref{A10}) (recall, also, (\ref{A9b})) it is evident that
for sufficiently large $|\lambda|$ the $n$ zeros of $P(r; \lambda)$ are distinct (i.e. simple). Therefore,
\begin{equation}
D_P(\lambda) \not\equiv 0.
\label{A10aa}
\end{equation}
As a small illustration, let us mention that in the case $n = 2$ we have $D_P(\lambda) = [u_1 (b; \lambda) + u_2'(b; \lambda)]^2 - 4$
(recall that in the Hill operator theory \cite{M-W} it is the quantity $u_1 (b; \lambda) + u_2'(b; \lambda)$ which is called discriminant).
Furthermore, in the Hill case the zeros of $D_P(\lambda)$ are simple or double and correspond to the periodic and antiperiodic eigenvalues
of $\mathcal{L}$. This correspondence is consistent in the sense that if $\lambda^\ast$ is a double zero of $D_P(\lambda)$, then the equation
$\mathcal{L} u = \lambda^\ast u$ has two linearly independent periodic or antiperiodic solutions (in some sense, we can say that algebraic
and geometric multiplicities are equal).

An additional observation is that, since $P(r; \lambda)$ is a monic polynomial, the discriminant of $P$ is \cite{L}
\begin{equation}
D_P(\lambda)
= \prod_{1 \leq j < k \leq n} \left[r_k(\lambda) - r_j(\lambda)\right]^2.
\label{A11}
\end{equation}
Thus, formula (\ref{A10}) implies that the entire function $D_P(\lambda)$ is of order $1 / n$. In view of (\ref{A9e}) and (\ref{A9g}), Corollary 2
applied to (\ref{A11}) yields the further corollary:

\smallskip

\textbf{Corollary 3.} Given a $\delta > 0$ there is a $K > 0$ such that if a zero $\lambda^{\ast}$ of $D_P(\lambda)$ satisfies
$|\lambda^{\ast}| > K$, then $(\lambda^{\ast})^{1/n} \in E_{\delta}$, where $E_{\delta}$ is the set defined in (\ref{A9g}).

\smallskip

The proof of Corollary 3 follows by observing that if $|\lambda^{\ast}|$ is large enough and $(\lambda^{\ast})^{1/n} \not\in E_{\delta}$,
then Corollary 2 and (\ref{A9e}) applied to (\ref{A11}) imply that $D_P(\lambda^{\ast})$ cannot vanish.

Now, let $\mathcal{M}(\mathbb{C})$ be the field of meromorphic functions on $\mathbb{C}$ (i.e. the quotient field of the integral
domain of entire functions). Then, $P(r; \lambda)$ is a polynomial over $\mathcal{M}(\mathbb{C})$; in other words,
$P$ belongs to the polynomial ring $\mathcal{M}(\mathbb{C})[r]$ (which is, of course, a unique factorization domain).
The fact that the discriminant $D_P(\lambda)$ of $P(r; \lambda)$ is not identically $0$ implies that the factorization
of $P(r; \lambda)$ over $\mathcal{M}(\mathbb{C})$ (as a polynomial in $r$) cannot have repeated factors. In fact, much more is true:

\smallskip

\textbf{Theorem 2.} The polynomial $P(r; \lambda)$ is irreducible over $\mathcal{M}(\mathbb{C})$.

\smallskip

\textit{Proof}. Let $r_{\ast}(\lambda)$ be a zero of $P(r; \lambda)$.
We will follow the analytic continuation of $r_{\ast}(\lambda)$ along a circular path
\begin{equation}
\lambda = (R/b)^n e^{i\theta}
\qquad\qquad
(\text{so that }\; R = b \left| \lambda^{1/n} \right|),
\label{A10c}
\end{equation}
where $R > 0$ is very large and $\theta$ is a real parameter.
Since the coefficients of $P(r; \lambda)$ are entire in $\lambda$ (and hence they do not change by analytic continuations),
it follows that any analytic continuation of a zero of $P(r; \lambda)$
will always remain a (possibly different) zero of $P(r; \lambda)$. In other words, analytic continuations of the zeros of $P(r; \lambda)$
can be viewed as elements of the Galois group $G(P / \mathcal{M}(\mathbb{C}))$ of $P(r; \lambda)$ over $\mathcal{M}(\mathbb{C})$.

When $\theta = 0$ we have that $r_{\ast}(\lambda)$ satisfies
(\ref{A10}); in particular, its leading behavior for large $R$ will be given by $\exp(\rho^j \lambda^{1 / n} b)$,
for some $j$. Now, we start moving $\theta$ away from $0$, in a continuous fashion, on the interval $[0, 2\pi]$.
Then, for each $\theta$ the quantity $r_{\ast}(\lambda)$, being a zero of $P(r; \lambda)$, will continue satisfying (\ref{A10}),
but with $\exp(\rho^k \lambda^{1 / n} b)$, where $k$ may be different from the $j$ which corresponded
to $\theta = 0$ (i.e. $k$ may depend on $\theta$). We claim that it is not possible to have two values $\theta_1$ and $\theta_2$ of $\theta$, with arbitrarily small distance
$|\theta_1 - \theta_2|$, such that (by invoking (\ref{A10})),
\begin{equation}
\left| r_{\ast}(\lambda_1) - \exp\left(\rho^j \lambda_1^{1 / n} b \right)\right|
\leq K \frac{\left|\exp\left(\rho^j \lambda_1^{1 / n} b \right) \right|}{R}
\label{A10a}
\end{equation}
and
\begin{equation}
\left| r_{\ast}(\lambda_2) - \exp\left(\rho^k \lambda_2^{1 / n} b \right)\right|
\leq K \frac{\left|\exp\left(\rho^k \lambda_2^{1 / n} b \right) \right|}{R}
\label{A10b}
\end{equation}
where $j \neq k$, $j,k \in \{0, 1, \dots, n-1\}$ and
\begin{equation*}
\lambda_1 := (R/b)^n e^{i \theta_1},
\qquad
\lambda_2 := (R/b)^n e^{i \theta_2}.
\end{equation*}
We will prove the claim by contradiction. Suppose (\ref{A10a}) and (\ref{A10b}) are satisfied for arbitrarily small
$|\theta_1 - \theta_2|$, namely for arbitrarily small $|\lambda_1 - \lambda_2|$. Since $r_{\ast}(\lambda)$ is analytic (hence continuous),
we can assume that the quantity $| r_{\ast}(\lambda_1) - r_{\ast}(\lambda_2) |$ is as close to $0$ as we wish.
%(e.g., $O(R^{-2})$).
Hence, if (\ref{A10a}) and (\ref{A10b}) are satisfied for arbitrarily small $|\theta_1 - \theta_2|$, we must have
(by taking $\lambda_2$ arbitrarily close to $\lambda_1$) that there is a $K > 0$ such that
\begin{equation}
\left| \exp\left(\rho^k \lambda_1^{1 / n} b\right) - \, \exp\left(\rho^j \lambda_1^{1 / n} b\right)  \right| \leq K \,
\frac{\max\left\{\left|\exp\left(\rho^k \lambda_1^{1 / n} b\right)\right|, \, \left| \exp\left(\rho^j \lambda_1^{1 / n} b\right)\right|\right\}}{R},
\label{A10bbb}
\end{equation}
where $R = |\lambda_1|^{1 / n} b$, while, as usual, $\lambda_1^{1 / n}$ is the principal branch of the $n$-th root and $\rho = e^{2\pi i / n}$.
However, by (\ref{A9e}), (\ref{A9g}), and Remark 2 we can see that there are arbitrarily large values of $R$ for which (\ref{A10bbb})
is not satisfied, and this is a contradiction.

%Without loss of generality we can assume that $|\exp(\rho^k \lambda_1^{1 / n} b)| \geq |\exp(\rho^j \lambda_1^{1 / n} b)|$
%(equivalently, $\cos(n^{-1}\theta_1 + 2\pi n^{-1} k) \geq \cos(n^{-1}\theta_1 + 2\pi n^{-1} j)$). Then,
%formula (\ref{A10bbb}) can be written as
%\begin{equation}
%\left| \exp\left((\rho^k - \rho^j) \lambda_1^{1 / n} b\right) - 1  \right| \leq \frac{K}{R}.
%\label{A10bb}
%\end{equation}
%If we set $z := (\rho^k - \rho^j) \lambda_1^{1 / n} b$ (so that $|z| = |\rho^k - \rho^j| R = c R$, where $c > 0$), we can see by using
%Lemma A1 of the Appendix that there are arbitrarily large values of $R$ for which (\ref{A10bb}) is not satisfied.
Therefore, we have reached the conclusion
that there are arbitrarily large $R$'s such that, as we analytically continue $r_{\ast}(\lambda)$
along the circular path $\lambda = (R/b)^n e^{i\theta}$, starting at $\theta = 0$, we always have that
\begin{equation}
\left| \frac{r_{\ast}(\lambda)}{\exp\left(\rho^j \lambda^{1 / n} b \right)} - 1 \right|
\leq \frac{K}{|\lambda|^{1 / n}}.
\label{A10c}
\end{equation}
Therefore, when $\theta$ reaches the value $2 \pi$, we arrive at a different zero $r_{\#}(\lambda)$ of $P(r; \lambda)$ satisfying
(at $\theta = 2 \pi$)
\begin{equation}
\left| \frac{r_{\#}(\lambda)}{\exp\left(\rho^{j+1} \lambda^{1 / n} b \right)} - 1 \right|
\leq \frac{K}{|\lambda|^{1 / n}}
\label{A10d}
\end{equation}
(as usual, $\lambda^{1 / n}$ denotes the principal branch of the $n$-th root). In particular, at $\theta = 2 \pi$
the leading behavior of $r_{\#}(\lambda)$ for large $R$ is $\exp(\rho^{j+1} \lambda^{1 / n} b)$.
If we continue the process of analytic continuation, by moving $\theta$, say, from $2 \pi$ to $4 \pi$,
we arrive at a new zero $r_{\bullet}(\lambda)$ of $P(r; \lambda)$ satisfying (at $\theta = 4 \pi$)
\begin{equation}
\left| \frac{r_{\bullet}(\lambda)}{\exp\left(\rho^{j+2} \lambda^{1 / n} b \right)} - 1 \right|
\leq \frac{K}{|\lambda|^{1 / n}}
\label{A10e}
\end{equation}
(at $\theta = 4 \pi$ the leading behavior of $r_{\bullet}(\lambda)$ for large $R$ is $\exp(\rho^{j+2} \lambda^{1 / n} b)$.
We can keep continuing analytically the zero of $P(r; \lambda)$ and we will always obtain a new zero until the $n$-th continuation where
we will obtain again $r_{\ast}(\lambda)$ (the zero we have started with). A crucial implication is that the Galois group
$G(P / \mathcal{M}(\mathbb{C}))$ is transitive, which is equivalent to the irreducibility of $P(r; \lambda)$ over $\mathcal{M}(\mathbb{C})$.
\hfill $\blacksquare$

\smallskip

Let us close this subsection by mentioning an interesting property of the Floquet multipliers, which has been proved in \cite{M}:

Suppose that $\mathcal{L}$ of (\ref{A1}) has real coefficients. Then $r_j(\lambda)$ is a Floquet multiplier of (\ref{A3}) if and only if
$r_j(\lambda)^{-1}$ is a Floquet multiplier of $\mathcal{L}^{\ast} v = \lambda v$. In other words, if $P(r; \lambda)$ is the characteristic
polynomial associated to $\mathcal{L}$ (see (\ref{A7})), then the characteristic polynomial $Q(r; \lambda)$ associated to the adjoint
$\mathcal{L}^{\ast}$ of the operator $\mathcal{L}$ is
\begin{equation}
Q(r; \lambda) = (-1)^n r^n P(1/r \, ; \lambda).
\label{A16}
\end{equation}
An immediate consequence of the above property is that in the case of the (real) self-adjoint operator $\mathcal{A}$ of (\ref{A1aa}),
the Floquet multipliers come in pairs of inverses. Equivalently, the associated characteristic polynomial $P(r; \lambda)$ is ``palindromic"
in the sense that it satisfies the equation
\begin{equation}
P(r; \lambda) = r^{2 \nu} P(1/r \, ; \lambda)
\qquad
(n = 2 \nu),
\label{A16ab}
\end{equation}
namely it has the form
\begin{equation}
P(r; \lambda)
= r^{2\nu} + A_{2\nu-1}(\lambda) r^{2\nu-1} + \cdots + A_1(\lambda) r + 1,
\label{A16aaa}
\end{equation}
where
\begin{equation}
A_{2\nu-l}(\lambda) = A_l(\lambda)
\qquad \text{for }\; l = 1, \dots, \nu-1.
\label{A17aaa}
\end{equation}
One can give a short proof of the above fact by using Theorem 2 and an idea of P. Kuchment \cite{K2} (see Remark 5 of Subsection 2.3).

\subsection{The Riemann surface $\Gamma$}

It is often more convenient (and natural), especially since $P(r; \lambda)$ is irreducible over $\mathcal{M}(\mathbb{C})$ (Theorem 2),
to view the Floquet multipliers $r_j(\lambda)$, $j = 1, \dots, n$, as the branches of a multi-valued
analytic function $r(\lambda)$, called the multiplier function. To make $r(\lambda)$ single-valued we need to introduce its
$n$-sheeted Riemann surface $\Gamma$ (on which $r$ is defined), where the sheets of $\Gamma$ are copies of the complex plane,
i.e. they do not contain $\infty$. Notice that the property of the Floquet multipliers mentioned at the end of the previous subsection implies
that, in the case of real coefficients, if $r(\lambda)$ is the multiplier of $\mathcal{L}$ (i.e. of (\ref{A3})) then the multiplier of
$\mathcal{L}^{\ast}$ is $r(\lambda)^{-1}$, and hence it ``lives" on the same Riemann surface $\Gamma$.

The surface $\Gamma$ can be identified with the complex (transcendental) curve over the $\lambda$-plane
(sometimes called the \textit{multiplier curve} of the operator $\mathcal{L}$)
\begin{equation}
P(r; \lambda) = 0.
\label{A14}
\end{equation}
In fact, one can use equations like (\ref{A3}) in order to obtain Riemann surfaces like $\Gamma$ whose genuses are, generically, infinite \cite{M},
\cite{Mc}, \cite{M-Mo}, \cite{M-T}.

For example, in the case Hill case ($n = 2$) we have that $\Gamma$ is the (generically infinite-genus)
hyperelliptic-type surface associated to the double-valued function
$\sqrt{D_P(\lambda)} = \sqrt{u_1 (b; \lambda) + u_2'(b; \lambda)]^2 - 4}$.
It follows that in the Hill case the projections on $\mathbb{C}$ of the ramification points of $\Gamma$ are the simple zeros of
$D_P(\lambda)$ and that $\Gamma$ can be compactified by adding $\infty$ to become a finite-genus compact Riemann surface $\overline{\Gamma}$
if and only if all but finitely many zeros of $D_P(\lambda)$ are double (this is the periodic finite-zone potential case).

For $n = 3$ the structure of $\Gamma$ has been studied in \cite{Mc}.
In the case of the ($4$-th order) periodic Euler-Bernoulli (or beam) operator the structure of $\Gamma$ has been studied in \cite{P1},
\cite{P2}, and \cite{P3}, while in the case of a general operator $\mathcal{L}$ the asymptotic structure (i.e. for $|\lambda|$ sufficiently large)
of $\Gamma$ has been determined in \cite{M0} and \cite{M}.

\smallskip

\textbf{Remark 3.} In the very special case where all coefficients of $\mathcal{L}$ are constant, say $p_k(x) \equiv \gamma_k$,
$k = 0, \dots, n-2$ (recall (\ref{A1}), (\ref{A0})), we have that $r(\lambda) = e^{w(\lambda) b}$, where $w = w(\lambda)$ satisfies
\begin{equation}
w^n + \sum_{k=0}^{n-2} \gamma_k w^k = \lambda.
\label{AN1}
\end{equation}
It follows that the compactification $\overline{\Gamma}$ of $\Gamma$ (by adding $\infty$) is the graph of the $n$-th degree
polynomial $\lambda = \lambda(w)$ appearing in (\ref{AN1}) and, therefore, $\overline{\Gamma}$
has genus $0$, hence, it is topologically and, furthermore, conformally equivalent to the sphere $\overline{\mathbb{C}}$ \cite{S1}.
For example, in the unperturbed case the obtained compactified surface is the
Riemann surface of $\lambda^{1 / n}$, hence it has only two $n$-th root ramification points, one at $0$ and the other at $\infty$,
while the multi-valued multiplier function is
\begin{equation}
\tilde{r}(\lambda) = \exp\left(\lambda^{1 / n} b\right),
\label{A14a}
\end{equation}
where, here $\lambda^{1 / n}$ denotes the multi-valued analytic $n$-th root function, and not just its principal branch, as before.
\hfill $\diamondsuit$

\smallskip

\textbf{Proposition 1.} Suppose that the operator $\mathcal{A}$ of (\ref{A1aa}) has constant coefficients.
Then the spectrum of $(-1)^{\nu}\mathcal{A}$ has the form $\sigma((-1)^{\nu} \mathcal{A}) = [\lambda_0, \infty)$,
where $\lambda_0$ is some real number.

\smallskip

\textit{Proof}. Consider the equation $(-1)^{\nu} \mathcal{A} u = \lambda u$ and assume $a_k(x) \equiv (-1)^{\nu} \alpha_k$ for
$k = 0, \dots, \nu - 1$. The Floquet multipliers associated to the aforementioned equation satisfy $r_j(\lambda) = e^{w_j(\lambda) b}$, where
$w_j(\lambda)$, $j = 1, \dots, 2 \nu$ are the roots of the polynomial equation
\begin{equation}
(-1)^{\nu} w^{2 \nu} + \sum_{k=0}^{\nu-1} \alpha_k w^{2k} = \lambda.
\label{AN1a}
\end{equation}
Then, (\ref{A13}) and self-adjointness imply that
\begin{equation*}
\sigma(\mathcal{A}) = \left\{ \lambda \in \mathbb{R} : \, \Re\{w_j(\lambda)\} = 0, \ \text{for some } j \right\}
= \left\{ \lambda \in \mathbb{R} : \, w_j^2(\lambda) \leq 0, \ \text{for some } j \right\}.
%\label{AN2}
\end{equation*}
By substituting $w^2 = \eta$ (\ref{AN1a}) becomes
\begin{equation*}
q(\eta) := (-1)^\nu \eta^\nu + \sum_{k=0}^{\nu-1} \alpha_k \eta^k = \lambda.
%\label{AN3}
\end{equation*}
Since $q(\eta) \to \infty$ as $\eta \to -\infty$ we get that
\begin{equation*}
\lambda_0 := \min_{\eta \leq 0} q(\eta)
\end{equation*}
exists in $\mathbb{R}$. It follows that, if $\lambda \geq \lambda_0$ there is a $\eta \leq 0$ such that $q(\eta) = \lambda$, while
for $\lambda < \lambda_0$ there is no such a $\eta$.
\hfill $\blacksquare$

\smallskip

In Remark 9 of Section 4 there is a discussion of how to formulate the converse of Proposition 1.

We continue with a definition of some special sets of points of the complex plane.

\smallskip

\textbf{Definition 1.} We denote by $\mathcal{R}_{\Gamma}$ the set of $\lambda \in \mathbb{C}$ over which lie the ramification points of $\Gamma$.
Also, we denote by $\mathcal{Z}$ the set of zeros of the discriminant $D_P(\lambda)$ of $P(r; \lambda)$
($\mathcal{Z}$ is countably infinite by the fact that $D_P(\lambda)$ has order $1/n$ and the Hadamard Factorization Theorem).
Finally, we denote by $\mathcal{Z}_J$ the set of $\lambda \in \mathbb{C}$ such that $T(\lambda)$ has a Jordan anomaly, i.e. it is not diagonalizable.

\smallskip

Of course, the set $\mathcal{R}_{\Gamma}$ is the set of branch points of the multiplier function $r(\lambda)$, since each branch point of
$r(\lambda)$ is the projection of some ramification point of $\Gamma$.
Notice that $\mathcal{R}_{\Gamma} \subset \mathcal{Z}$, while ``generically" we have
$\mathcal{R}_{\Gamma} = \mathcal{Z}$. However, in general $\mathcal{Z}$ and $\mathcal{R}_{\Gamma}$ may not be equal, since,
it can happen that all points of $\Gamma$ lying over a specific number $\lambda^{\star} \in \mathcal{Z}$, are unramified. For example
in the unperturbed case we have (see (\ref{A9d}))
\begin{equation*}
\tilde{\mathcal{Z}} = \{ \tilde{\lambda}_{m,d} :\;  m \in \mathbb{N}, \; d = 1, \dots, n-1 \} \cup \{0\},
\quad \text{while} \quad
\mathcal{R}_{\tilde{\Gamma}} = \{ 0 \}.
\end{equation*}
Also, $\mathcal{Z}_J \subset \mathcal{Z}$ and, again, the two sets may not be equal (e.g., in the case of $\tilde{\mathcal{L}}$, where
$\tilde{\mathcal{Z}}_J = \{ 0 \}$, as we will see below, in the end of the next subsection), although, ``generically" we have
$\mathcal{Z}_J = \mathcal{Z}$. Finally, if Conjecture 1 of Section 4 is true, then $\mathcal{R}_{\Gamma} \subset \mathcal{Z}_J$.

The behavior of the large (in absolute value) elements of the set $\mathcal{Z}$ is described in Corollary 3.

\smallskip

\textbf{Remark 4.} The multiplier $r(\lambda)$ is a special Baker-Akhiezer-type function (see, e.g., \cite{B-K}) defined on $\Gamma$.
It never vanishes and its branches can be chosen to satisfy (\ref{A10}). We can think of it as an analog of
the exponential function $e^{c \lambda}$ (of $\mathbb{C}$) for $\Gamma$.

The derivative $r'(\lambda)$ (here, the prime denotes derivative with respect to $\lambda$ since there is no risk of confusion) is meromorphic on $\Gamma$ and its poles are located at the ramification points of $\Gamma$. More specifically, by using a local uniformizing parameter we can
easily see that $\ell = (\lambda, r) \in \Gamma$ is a pole of $r'$ of order $p$ if and only if
$\ell$ is a ramification point of $\Gamma$ of index $p$ (recall that the index of a ramification point $\ell$ is $1$ less than
its degree, i.e. $1$ less than the number of adjacent sheets with common point $\ell$). The asymptotics of $r_j'(\lambda)$, $j = 1, \dots, n$,
as $\lambda$ approaches $\infty$ staying away from the branch points of $r(\lambda)$ (i.e. staying away from the set $\mathcal{R}_{\Gamma}$), can
be computed by Cauchy's formula applied to (\ref{A10}):
\begin{equation}
r_j'(\lambda) = \frac{1}{2 \pi i} \oint_C \frac{r_j(z)}{(z - \lambda)^2} dz
= \frac{1}{2 \pi i} \oint_C \frac{exp\left(\omega_k z^{1/n} b\right)}{(z - \lambda)^2}
\left[1 + O\left(\lambda^{-1 / n}\right)\right] dz,
\label{A15aaa}
\end{equation}
where $C$ is a circle of a fixed radius $a > 0$, centered at $\lambda$ and not enclosing $0$, nor any branch points of $r(\lambda)$, while
$\omega_k$ is an $n$-th root of $1$ as in Corollary 2. From (\ref{A15aaa}) we get,
\begin{equation}
r_j'(\lambda) = \frac{\omega_k b}{n \lambda^{(n-1)/n}}\exp\left(\omega_k \lambda^{1/n} b\right) \left[1 + O\left(\lambda^{-1 / n}\right) \right]
\label{A15aa}
\end{equation}
as $\lambda$ approaches $\infty$ staying one unit away from all points of the set $\mathcal{R}_{\Gamma}$.

The function $(r' / r)(\lambda)$ is also meromorphic on $\Gamma$ and its poles are of the same nature as the poles of $r'(\lambda)$, since
$r(\lambda)$ does not vanish on $\Gamma$. Furthermore (\ref{A10}) and (\ref{A15aa}) imply
\begin{equation}
\frac{r_j'(\lambda)}{r_j(\lambda)} = \frac{\omega_k b}{n \lambda^{(n-1)/n}}
\left[1 + O\left(\lambda^{-1 / n}\right) \right],
\qquad
j = 1, \dots, n,
\label{A15}
\end{equation}
as $\lambda$ approaches $\infty$ staying one unit away from all points of the set $\mathcal{R}_{\Gamma}$
(as in Corollary 2, if $\lambda^{1/n} \not\in E_{\delta}$, to each branch $r_j' /r_j$ corresponds a different $\omega_k$).
If $\Gamma$ has finitely many ramification points, then (\ref{A15}) is true as $\lambda$ approaches $\infty$ without any restrictions. Therefore,
(\ref{A15}) implies that $(r' / r)(\lambda)$ extends meromorphically
to the compactification $\overline{\Gamma} = \Gamma \cup \{\infty\}$ (notice that by (\ref{A15}) we have $(r' / r)(\infty) = 0$).
Thus, in this case $\overline{\Gamma}$ is a compact
Riemann surface (and, hence, the genus $g$ of $\overline{\Gamma}$ is finite), $\infty$ is a $n$-th root branch point,
and $(r' / r)(\lambda)$  is a (purely) algebraic function
(in other words, there is a polynomial $Q(\cdot\, , \cdot)$ in two variables such that $Q(r' / r, \lambda) = 0$). As we have seen,
$(r' / r)(\lambda)$  is a very special meromorphic function on $\overline{\Gamma}$. Its poles are completely determined by
the ramification points of $\overline{\Gamma}$, while $\infty$ is a zero of $r' / r$ of multiplicity $n-1$.
In that case the Riemann-Hurwitz formula (see, e.g., \cite{M-M} or \cite{S1}) applies to $\overline{\Gamma}$.
The formula reads
\begin{equation}
R_{\text{tot}} = 2(n + g -1),
\label{A15a}
\end{equation}
where $R_{\text{tot}} $ is the total ramification index, i.e. the sum of the indices of all ramifications.
Notice that by (\ref{A15}) the index of the ramification at $\infty$ is $n - 1$.
For example, since $g \geq 0$, formula (\ref{A15a}) implies that $R_{\text{tot}}  \geq 2(n - 1)$, which, in particular, tells us
that $\Gamma$ always must have some ramification point(s).
\hfill $\diamondsuit$

\smallskip

A well-known example where $\Gamma$ has finitely many ramification points, hence $\overline{\Gamma}$ is a surface with a genus $g$,
is the case where $\mathcal{L}$ is a Hill operator (thus $n = 2$) with a
finite-zone, $b\,$-periodic (real) potential. Here it is known that the sets $\mathcal{R}_{\Gamma}$ and $\mathcal{Z}_J$ are equal
(in fact, this is always true in the Hill case) and,
by the Riemann-Hurwitz formula their cardinality is $2g + 1$ (by the definition of $\mathcal{R}_{\Gamma}$ and $\mathcal{Z}_J$, $\infty$ is
excluded from being an element of these sets).

\subsection{Floquet solutions}

Let $\lambda \in \mathbb{C}$ be fixed. Then, to each branch $r_j = r_j(\lambda)$, $j = 1, 2, \dots, n$, of the multiplier $r$ corresponds
at least one (up to linear independence) eigenvector $f_j$ of $T$, namely a solution of (\ref{A3}) such that
\begin{equation}
(\mathcal{T} f_j)(x; \lambda) = f_j(x+b; \lambda) = r_j(\lambda) f_j(x; \lambda)
\qquad \text{for all }\; x \in \mathbb{R}.
\label{A8}
\end{equation}
These eigenvectors are called (pure) Floquet solutions and it is not hard
to see that they have the form
\begin{equation}
f_j(x) = e^{i \kappa_j x} p_j(x),
\qquad \text{where }\; e^{i \kappa_j b} = r_j
\quad \text{and} \quad
p_j(x+b) = p_j(x),
\label{A9}
\end{equation}
where the $i \kappa_j$'s are the \textit{Floquet exponents} (the complex quantity $\kappa_j$ is sometimes called
\textit{rotation number} for rather obvious reasons, or
\textit{quasimomentum}, since in the case where $\mathcal{L}$ is a Schr\"{o}dinger operator it has a quantum mechanical meaning).

If $T(\lambda)$ is diagonalizable, then there are $n$ linearly independent
Floquet solutions. However, if for some $\lambda = \lambda^{\star}$ the
polynomial $P(r; \lambda^{\star})$ has multiple zeros, then
$T(\lambda^{\star})$ may not be diagonalizable, in other words it may have
Jordan anomalies (i.e. we may have $\lambda^{\star} \in \mathcal{Z}_J$---see Definition 1), hence there would be less than
$n$ linearly independent Floquet solutions. To be more precise, suppose that, for some $\lambda^{\star}$,
the matrix $T(\lambda^{\star})$ has $m$ distinct eigenvalues
$r_1, \dots, r_m$, with $m < n$. Let $f_{1j}, \dots, f_{n_jj}$ be the
eigenvectors (i.e. the pure Floquet solutions) that correspond to $r_j$,
$j = 1, \dots, m$. If $n_j > 1$ we say that we have \textit{coexistence} of $n_j$
(pure) linearly independent Floquet solutions associated to the multiplier $r_j$. Then, in the Jordan canonical
form of $T(\lambda^{\star})$ there are $n_j$ Jordan blocks whose diagonal
entries are all equal to $r_j$. Assume that the sizes of these Jordan blocks are
$s_{1j} \times s_{1j}, \dots, s_{n_jj} \times s_{n_jj}$ (where the $i$-th block
corresponds to the eigenvector $f_{ij}$). Then, to each $f_{ij}$
there correspond $s_{ij} - 1$ generalized eigenvectors
$g^l_{ij}$, $l = 2, \dots, s_{ij}$ ($l$ is a superscript),  such that
\begin{equation}
(\mathcal{T} g^l_{ij})(x) = g^l_{ij}(x+b)
= r_j g^l_{ij}(x) + g^{l-1}_{ij}(x),
\quad \text{with} \quad
g^1_{ij}(x) := f_{ij}(x).
\label{A9aa}
\end{equation}
Notice that, as an eigenvalue of $T(\lambda^{\star})$, the geometric multiplicity of $r_j(\lambda^{\star})$ equals $n_j$, while the algebraic
multiplicity of $r_j(\lambda^{\star})$ is $s_{1j} + \cdots + s_{n_jj}$. If $n_j > 1$, but $s_{ij} > 1$ for some $i$ (so that the algebraic
multiplicity of $r_j(\lambda^{\star})$ is greater than $n_j$), we can say that we
have \textit{partial coexistence} of (pure) linearly independent Floquet solutions associated to $r_j(\lambda^{\star})$,
while if $s_{ij} = 1$ for all $i = 1, \dots, n_j$, then the algebraic and geometric multiplicities of $r_j(\lambda^{\star})$ coincide and
we can say that we have \textit{complete} or \textit{full coexistence} of (pure) linearly independent Floquet solutions associated to
$r_j(\lambda^{\star})$. Clearly,
\begin{equation}
\sum_{j=1}^m \sum_{i=1}^{n_j} s_{ij} = n
\label{A9aaa}
\end{equation}
and $T(\lambda)$ is diagonalizable (for a fixed $\lambda \in \mathbb{C}$) if and only if $s_{ij} = 1$ for all $i, j$
(i.e. we have complete coexistence of pure Floquet solutions associated to $r_j(\lambda)$, for every $j = 1, 2, \dots, m$).

The $g^l_{ij}(x)$'s ($l > 1$) of (\ref{A9aa}) can be called generalized Floquet solutions of rank $l$ (the rank $l=1$ signifies
pure Floquet solution). Linear combinations of $g^l_{ij}(x)$'s corresponding to the same $j$ and $l$ could be, also, called
generalized Floquet solutions of rank $l$. Let us make this terminology official via the following (recursive) definition.

\smallskip

\textbf{Definition 2.} Consider the equation (\ref{A3}) for a fixed $\lambda \in \mathbb{C}$.
 A (pure) Floquet solution $u^1(x)$, namely a nontrivial solution satisfying
\begin{equation*}
\mathcal{(T} u^1)(x) = u^1(x+b)= r u^1(x),
\end{equation*}
is called a \textit{Floquet solution of rank} $1$ \textit{associated to the multiplier} $r$.
For $l \geq 2$, a solution $u^l(x)$ of (\ref{A3}) is called a (\textit{generalized})
\textit{Floquet solution of rank} $l$ \textit{associated to the multiplier} $r$, if there is a Floquet solution $u^{l-1}(x)$ of (\ref{A3})
of rank $(l-1)$ associated to $r$ such that
\begin{equation*}
(\mathcal{T} u^l)(x) = u^l(x+b)= r u^l(x) + u^{l-1}(x).
\end{equation*}
For completeness, the trivial solution $u^0(x) \equiv 0$ is called the \textit{Floquet solution of rank} $0$.

\smallskip

Of course, the number of linearly independent Floquet solution of rank $l \geq 2$ associated to a multiplier $r$ cannot exceed
the number of linearly independent Floquet solution of rank $(l-1)$ associated to $r$.

\smallskip

\textbf{Example 1.} For some $\lambda^{\star} \in \mathbb{C}$ let us assume that the Floquet matrix $T(\lambda^{\star})$ has $m < n$ distinct
eigenvalues $r_1(\lambda^{\star}), \dots, r_m(\lambda^{\star})$. Suppose that $r_1(\lambda^{\star})$ is a zero of $P(r; \lambda^{\star})$
of multiplicity $8$ (i.e. the algebraic multiplicity of $r_1(\lambda^{\star})$ is $8$), and in the Jordan canonical form of $T(\lambda^{\star})$
there are four Jordan blocks $J_1, \dots, J_4$ corresponding to $r_1 = r_1(\lambda^{\star})$:
\begin{equation*}
J_1 = [r_1],
\qquad
J_2 = \left[
\begin{array}{cc}
r_1 & 1 \\
0 & r_1
\end{array}
\right],
\qquad
J_3 = \left[
\begin{array}{cc}
r_1 & 1 \\
0 & r_1
\end{array}
\right],
\qquad
J_4 = \left[
\begin{array}{ccc}
r_1 & 1 & 0  \\
0 & r_1 & 1 \\
0 & 0 & r_1
\end{array}
\right]
\end{equation*}
(hence the geometric multiplicity of $r_1(\lambda^{\star})$ is $4$, i.e. $n_1 = 4$, and also $s_{11} = 1$, $s_{21} = 2$,
$s_{31} = 2$, and $s_{41} = 3$). Then, the space of Floquet solutions of rank $1$
associated to $r_1$ has dimension $4$ (e.g., $f_{11}, f_{21}, f_{31}, f_{41}$ is a basis for this space),
while the space of Floquet solutions of rank $2$ associated to $r_1$ has
dimension $3$ (e.g., $g^2_{21}, g^2_{31}, g^2_{41}$ is a basis for this space) and the space of Floquet solutions
of rank $3$ associated to $r_1$ has dimension $1$ (and it is generated by $g^3_{41}$). There are no Floquet solutions of rank
$l \geq 4$ associated to $r_1$. We have partial coexistence since there are four pure Floquet solutions corresponding
to $r_1(\lambda^{\star})$, while the algebraic multiplicity of to $r_1(\lambda^{\star})$ is $8$.

\smallskip

\textbf{Remark 5.} As we have mentioned at the end of Subsection 2.1, if $\mathcal{A}$ is a (real) self-adjoint operator (as in (\ref{A1aa})),
then its associate characteristic polynomial $P(r; \lambda)$ is ``palindromic", i.e. it satisfies formula (\ref{A16ab}). A short proof of this fact
(based on Theorem 2 and an idea of P. Kuchment \cite{K2}), which exploits the $\lambda$-dependence goes as follows.

Let $\lambda \in \sigma(\mathcal{A}) \setminus \mathcal{Z}$ (see Definition 1). Then, by (\ref{A13}) there is a Floquet solution $f_j(x)$ such that
\begin{equation}
\mathcal{A} f_j = \lambda f_j,
\qquad
f_j(x + b) = r_j(\lambda) f_j(x)
\label{FFT1}
\end{equation}
with $|r_j(\lambda)| = 1$. We exclude the countably many $\lambda$'s for which $r_j(\lambda) = \pm 1$ and we take complex conjugates in (\ref{FFT1})
to get (remember that the spectrum of $\mathcal{A}$ is real)
\begin{equation*}
\mathcal{A} \bar{f}_j = \lambda \bar{f}_j,
\qquad
\bar{f}_j(x + b) = \overline{r_j(\lambda)} \, \overline{f_j(x)} = r_j(\lambda)^{-1} \overline{f_j(x)}.
%\label{FFT2}
\end{equation*}
Hence, $\overline{f_j(x)}$ is also a (pure) Floquet solution with corresponding multiplier $r_j(\lambda)^{-1}$. Therefore, if
$\lambda \in \sigma(\mathbb{A}) \setminus \mathcal{Z}$ and $r_j(\lambda)$ is a Floquet multiplier, so is $r_j(\lambda)^{-1}$.
But then, analytic continuation and Theorem 2 (namely the irreducibility of $P(r; \lambda)$) imply that this holds for all $\lambda \in \mathbb{C}$.

It is remarkable that the original proof of the more general fact mentioned at the end of Subsection 2.1, regarding the real coefficient operators
$\mathcal{L}$ and $\mathcal{L}^{\ast}$, does not make use of the $\lambda$-dependence.
\hfill $\diamondsuit$

\smallskip

It is desirable to view the Floquet solutions $f_j(x; \lambda)$, $j = 1, \dots, n$ (for any $x$) as the branches
of a multi-valued analytic function. Now, each Floquet solution $f_j(x; \lambda)$ is determined up to
an $x$-independent factor. In order to overcome this ambiguity we introduce the normalized Floquet solutions \cite{Mc}, \cite{M}, \cite{No}
\begin{equation}
\phi_j(x; \lambda) := \frac{f_j(x; \lambda)}{f_j(0; \lambda)},
\label{A9a}
\end{equation}
so that $\phi_j(0; \lambda) = 1$. To be more specific, given a branch $r_j(\lambda)$ of the (multi-valued) multiplier $r(\lambda)$
let us consider the normalized Floquet solution $\phi_j(x; \lambda)$ of (\ref{A3}) such that
\begin{equation}
\phi_j(x + b; \lambda) = r_j(\lambda) \phi_j(x; \lambda),
\qquad
\phi_j(0; \lambda) = 1
\label{NF1}
\end{equation}
(the existence of $\phi_j(x; \lambda)$ will become clear in the sequel). We expand $\phi_j$ in terms of the fundamental solutions
\begin{equation}
\phi_j(x; \lambda) = u_1(x; \lambda) + c_{j2}(\lambda) u_2(x; \lambda) + \cdots + c_{jn}(\lambda) u_n(x; \lambda).
\label{NF2}
\end{equation}
Thus, the vector
\begin{equation}
\phi_j = \left[1, c_{j2}(\lambda), \dots, c_{jn}(\lambda) \right]^{\top}
\label{NF3}
\end{equation}
is an eigenvector of $T(\lambda)$ with eigenvalue $r_j(\lambda)$. Therefore (see (\ref{A6})) we can use Cramer's rule, to express
$c_{jk}$, $k = 2, 3, \dots, n$, in terms of $r_j(\lambda)$ and the $t_{ij}(\lambda)$'s. Let us denote by $T_{11}(\lambda)$ the
$(n-1) \times (n-1)$ matrix obtained from $T(\lambda)$ of (\ref{A6}) after erasing the first row and the first column, i.e.
\begin{equation}
T_{11} = T_{11}(\lambda)
= \left[t_{ij}(\lambda)\right]_{2 \leq i,j \leq n}.
\label{NF4}
\end{equation}
Next, we introduce the following $(n-1)$-th degree polynomial over the function field $\mathcal{M}(\mathbb{C})$:
\begin{equation}
\Delta(r) = \Delta(r; \lambda) := \det [T_{11}(\lambda) - r I],
\label{NF5}
\end{equation}
where $I$ is the $(n-1) \times (n-1)$ identity matrix. Notice that, if $r_j(\lambda)$ is a branch of the multiplier $r(\lambda)$,
then $\Delta(r_j(\lambda); \lambda) \not \equiv 0$, since by Theorem 2 $r_j(\lambda)$ has degree $n$ over the field $\mathcal{M}(\mathbb{C})$,
while $\deg \Delta(r) = n-1$.
For typographical convenience we denote by $[T_{11}(\lambda) - r I]_k$, $k = 2, 3, \dots, n$, the $(n-1) \times (n-1)$ matrix obtained
from $[T_{11}(\lambda) - r I]$ after replacing its $k$-th column by the column vector
\begin{equation*}
\left[-t_{21}(\lambda), -t_{31}(\lambda), \dots, -t_{n1}(\lambda) \right]^{\top}.
\end{equation*}
Then, Cramer's rule yields
\begin{equation}
c_{jk}(\lambda) =\frac{\det [T_{11}(\lambda) - r_j(\lambda) I]_k}{\Delta(r_j(\lambda))},
\qquad
k = 2, 3, \dots, n.
\label{NF6}
\end{equation}
It follows that for each $k = 2, \dots, n$, the quantities $c_{jk}(\lambda)$, $j = 1, \dots, n$, are the branches of a meromorphic function
on $\Gamma$ expressed as
\begin{equation}
c_{k}(\lambda) =\frac{\det [T_{11}(\lambda) - r(\lambda) I]_k}{\Delta(r(\lambda))},
\qquad
k = 2, 3, \dots, n,
\label{NF7}
\end{equation}
where $r(\lambda)$ is the (multi-valued) multiplier function. If $\ell = (\lambda, r) \in \Gamma$ is a pole of $c_{k}(\lambda)$
and $\lambda^{\star}$ is the projection of $\ell$ on $\mathbb{C}$, then $\Delta(r(\ell); \lambda^{\star}) = 0$ (however, in general the
converse may not be true). In other words, the projections of the poles of the $\phi_j(x; \lambda)$'s are zeros of the entire function
\begin{equation}
\prod_{j=1}^n \Delta(r_j(\lambda); \lambda).
\label{NF55}
\end{equation}

Finally, using (\ref{NF6}) in (\ref{NF2}) we can see that each $\phi_j(x; \lambda)$, $j = 1, 2, \dots, n$, can be seen as a branch of a
meromorphic normalized Floquet solution $\phi(x; \lambda)$ defined on $\Gamma$, such that $\phi(0; \lambda) = 1$, provided $\lambda$
is not a pole of $\phi(x; \lambda)$. These poles, being among the zeros of $\Delta(r(\lambda); \lambda)$, are independent of $x$.
Equation (\ref{NF2}) can be written in the form
\begin{equation}
\phi(x; \lambda) = u_1(x; \lambda) + c_{2}(\lambda) u_2(x; \lambda) + \cdots + c_{n}(\lambda) u_n(x; \lambda).
\label{NF10}
\end{equation}
Let us, also, mention that, as a function of $\lambda$, $\phi(x; \lambda)$ is a Baker-Akhiezer function (see, e.g., \cite{B-K}). In fact,
some authors \cite{Mc}, \cite{M} call $\phi(x; \lambda)$ \textit{the} Baker-Akhiezer function.

In the unperturbed case we have
\begin{equation}
\tilde{\phi}(x; \lambda) = \exp\left(\lambda^{1/n} x \right)
\label{NF9}
\end{equation}
(in particular, $\tilde{\phi}(x; 0) \equiv 1$), where $\lambda^{1 / n}$ denotes the multi-valued analytic $n$-th root function,
and not just its principal branch. The branches of $\tilde{\phi}(x; \lambda)$ are
\begin{equation}
\tilde{\phi}_j(x; \lambda) = \exp\left( \rho^{j-1} \lambda^{1/n} x \right),
\qquad
j = 1, 2, \dots, n,
\label{NF8}
\end{equation}
where, here $\lambda^{1 / n}$ denotes the principal branch of the $n$-th root function, while $\rho$ is given by (\ref{A6b}).
It follows that, in the unperturbed case (see (\ref{NF8})), if $\lambda \ne 0$, then $\tilde{T}(\lambda)$ has $n$ distinct eigenvectors and hence it
is diagonalizable. For the remaining value of $\lambda$, namely $\lambda = 0$, we get that $\tilde{r}_j(0) = 1$ for all $j$, $\tilde{T}(0)$
is similar to a $n \times n$ Jordan block and hence there is only one (pure)
Floquet solution, namely the constant function (hence, $\tilde{\mathcal{Z}}_J = \{0\}$---see Definition 1). The generalized Floquet solution of rank $l$, with $2 \leq l \leq n$, is a polynomial of degree $l-1$.

\subsection{The asymptotic behavior of $\phi(x; \lambda)$ as $\lambda \to \infty$}

The following important estimate regarding the behavior of $\phi(x; \lambda)$ as $\lambda \to \infty$ has been derived in \cite{M} via a quite
technical computation (in the special case $n = 3$ the estimate had been already derived in \cite{Mc} with the help of an idea of J. Moser).

Let $\delta > 0$ and, as usual, $\lambda = \zeta^n$, while $\omega_1, \dots, \omega_n$ are the $n$-th roots of $1$. Then
\begin{equation}
\phi(x; \lambda) = e^{\omega_k \zeta x} [1 + o(1)]
\qquad \text{as }\;
\zeta \to \infty, \quad \zeta \not\in E_{\delta},
\label{dSM}
\end{equation}
where $E_{\delta}$ is the set introduced in (\ref{A9g}) (in \cite{M} the coefficients of $\mathcal{L}$ are assumed real-valued;
however, this assumption is not used in the derivation of the estimate (\ref{dSM})).

\subsection{The ``shifted" operators}

For $\xi \in \mathbb{R}$ we consider the operator
\begin{equation}
\mathcal{L}_{\xi} u := \frac{d^n u}{dx^n} + \sum_{k=0}^{n-2} p_k(x + \xi) \frac{d^k u}{dx^k},
\label{A19}
\end{equation}
Of course, $\mathcal{L}_{\xi + b} = \mathcal{L}_{\xi}$ (with $\mathcal{L}_0 = \mathcal{L}$) and $v(x)$ satisfies $\mathcal{L}_{\xi} v = \lambda v$ if and only if $v(x) = u(x + \xi)$, where $u(x)$
satisfies $\mathcal{L} u = \lambda u$. Thus, Floquet solutions of $\mathcal{L}_{\xi} v = \lambda v$ are in one to one correspondence
with Floquet solutions of $\mathcal{L}u = \lambda u$ and the associated multipliers are equal. Furthermore, this is also true for generalized Floquet solutions. Therefore, the Floquet matrix $T_{\xi}(\lambda)$ of $\mathcal{L}_{\xi}$ is similar
to the Floquet matrix $T(\lambda)$ of $\mathcal{L}$. In particular, $T_{\xi}(\lambda)$ and $T(\lambda)$ have the same characteristic
polynomial $P(r; \lambda)$.

Now suppose that each coefficient $p_k(x)$ of $\mathcal{L}$, $k = 0, \dots, n-2$, extends to a function $p_k(z)$, $z = x +iy$, which is analytic
in a domain of the complex plane containing the real axis (for example, each $p_k(x)$ is a $b$-periodic trigonometric polynomial). Since
$p_k(x + b) = p_k(x)$, analytic continuation implies that $p_k(z)$ must be analytic in some open strip $S$ containing $\mathbb{R}$ and
$p_k(z + b) = p_k(z)$ for all $z \in S$. Then, $\mathcal{L}_{\xi}$ of (\ref{A19}) makes sense for any $\xi \in S$, has $b$-periodic coefficients,
and its associated fundamental solutions $u_j(x; \lambda; \xi)$, $j = 1, \dots, n$, are analytic in $\xi$. It follows that the coefficients
of the characteristic polynomial $P(r; \lambda; \xi)$ of the matrix $T_{\xi}(\lambda)$ are, also, analytic in $\xi$. However, as we saw the
coefficients of $P(r; \lambda; \xi)$ do not depend on $\xi$, if $\xi$ is real. Hence, by analytic continuation they do not depend on $\xi$ for
any $\xi \in S$. In other words, $P(r; \lambda; \xi)$ is independent of $\xi$ and, therefore, $\Gamma_{\xi} = \Gamma$ for any $\xi \in S$,
where by $\Gamma_{\xi}$ we denote the Riemann surface associated to $\mathcal{L}_{\xi}$. It follows that we have constructed easily a two (real)
parameter family of operators $\mathcal{L}_{\xi}$, $\xi \in S$, with a common associated Riemann surface $\Gamma$ and a common
polynomial $P(r; \lambda)$.

\smallskip

\textbf{Remark 6.} Formula (\ref{TA6}) holds for the operators $\mathcal{L}_{\xi}$, $\xi \in \mathbb{R}$, with a constant $K$ which is
independent of $\xi \in \mathbb{R}$ or, in the case of analytic coefficients, as long as $\Im\{\xi\}$ stays bounded.
\hfill $\diamondsuit$

\smallskip

Finally, let us mention that certain properties, which can be proved for the case of analytic coefficients, can be also
extended by density arguments to more general cases where the operator $\mathcal{L}$ has less smooth coefficients.

\section{A multipoint eigenvalue problem}

We introduce the following multipoint problem for the operator $\mathcal{L}$ as an analog
of the Dirichlet problem for the Hill's equation:
\begin{equation}
\mathcal{L} u = \lambda u,
\qquad
u(0) = u(b) = u(2b) = \cdots = u((n-1)b) = 0
\label{B1}
\end{equation}
This unorthodox boundary value problem first appeared in \cite{Mc} (its introduction there was suggested by L. Nirenberg) for the case
$n = 3$ and later, independently, in \cite{P1} for the case of the periodic Euler-Bernoulli operator (where $n = 4$).

An eigenvalue of (\ref{B1}) is a $\lambda \in \mathbb{C}$ for which (\ref{B1})
has a nontrivial solution. Any such solution is an eigenfunction of (\ref{B1}).
Since every solution of (\ref{B1}) is a linear combination of the fundamental
solutions $u_j(x; \lambda)$, $j = 1, \dots, n$, it follows that $\lambda$ is
an eigenvalue of (\ref{B1}) if and only if $\lambda$ is a zero of the entire
function
\begin{equation}
H(\lambda) := \left\vert
\begin{array}{cccc}
u_1(0; \lambda) & u_2(0; \lambda) & \cdots & u_n(0; \lambda)  \\
u_1(b; \lambda) & u_2(b; \lambda) & \cdots & u_n(b; \lambda) \\
\vdots & \vdots & \ddots & \vdots  \\
u_1((n-1)b; \lambda) & u_2((n-1)b; \lambda) & \cdots & u_n((n-1)b; \lambda)
\end{array}
\right\vert .
\label{B2}
\end{equation}
Notice that the order of $H(\lambda)$ is at most $1/n$. Since $u_1(0; \lambda) = 1$ and
$u_2(0; \lambda) = u_3(0; \lambda) = \cdots = u_n(0; \lambda) = 0$,
$H(\lambda)$ takes the form
\begin{equation}
H(\lambda) = \left\vert
\begin{array}{cccc}
u_2(b; \lambda) & u_3(b; \lambda) & \cdots & u_n(b; \lambda)  \\
u_2(2b; \lambda) & u_3(2b; \lambda) & \cdots & u_n(2b; \lambda) \\
\vdots & \vdots & \ddots & \vdots  \\
u_2((n-1)b; \lambda) & u_3((n-1)b; \lambda) & \cdots & u_{2\nu}((n-1)b; \lambda)
\end{array}
\right\vert .
\label{B3}
\end{equation}
Now, let $\mu \in \mathbb{C}$ be such that $H(\mu) = 0$. We denote by $m_a(\mu)$
the multiplicity of $\mu$ as a zero of $H$. We refer to $m_a(\mu)$ as the
\textit{algebraic multiplicity of} $\mu$. The eigenfunctions of (\ref{B1})
associated to $\lambda = \mu$ form a vector space ${\cal D}_{\mu}$ (the
eigenspace associated to $\mu$) with dimension
$\dim {\cal D}_{\mu} =: m_g(\mu)$, the \textit{geometric multiplicity of} $\mu$.
For the sake of completeness we can define these multiplicities for any $\mu \in \mathbb{C}$
by setting $m_a(\mu) = m_g(\mu) = 0$ whenever $\mu$ is not an eigenvalue of (\ref{B1}), i.e.
whenever $H(\mu) \ne 0$. Needless to say that $m_a(\mu) > 0$ if and only if $m_g(\mu) > 0$.

\smallskip

\textbf{Observation.} Suppose $\mu$ is an eigenvalue of (\ref{B1}) with geometric multiplicity $m_g(\mu)$. Then, the $n \times n$
matrix
\begin{equation}
M(\lambda) := \left[
\begin{array}{cccc}
1 & 0 & \cdots & 0  \\
u_1(b; \lambda) & u_2(b; \lambda) & \cdots & u_n(b; \lambda) \\
\vdots & \vdots & \ddots & \vdots  \\
u_1((n-1)b; \lambda) & u_2((n-1)b; \lambda) & \cdots & u_n((n-1)b; \lambda)
\end{array}
\right],
\label{B2a}
\end{equation}
whose determinant appears in (\ref{B2}), has rank $n - m_g(\mu)$. It follows that all derivatives of $H(\lambda)$ up to order
$m_g(\mu) - 1$ vanish at $\lambda = \mu$. Hence,
\begin{equation}
m_a(\mu) \geq m_g(\mu).
\label{B3a}
\end{equation}
In particular, if $m_a(\mu) = 1$, then $m_g(\mu) = 1$.

\smallskip

\textbf{Theorem 3.} For any eigenvalue $\mu$ of (\ref{B1}) its associated
eigenspace ${\cal D}_{\mu}$ has a basis consisting of pure and (possibly) generalized Floquet solutions. If a
generalized Floquet solution $\psi^l(x)$ of rank $l \geq 2$ associated to the multiplier $r_j(\mu)$ (see Definition 2) is in ${\cal D}_{\mu}$,
then the Floquet solution $\psi^{l-1}(x)$ of rank $(l-1)$ associated to $r_j(\mu)$ such that
\begin{equation*}
(\mathcal{T} \psi^l)(x) = \psi^l(x+b) = r_j(\mu) \psi^l(x) + \psi^{l-1}(x)
\end{equation*}
is in ${\cal D}_{\mu}$ too.

\smallskip

\textit{Proof}. Let $\mu$ be a zero of $H(\lambda)$ and suppose that the Floquet
matrix $T(\mu)$ has $m$ distinct eigenvalues $r_1, \dots, r_m$ ($m \leq n$). As
in Subsection 2.2 for each $j = 1, \dots, m$, we consider the (pure) Floquet
solutions $f_{1j}, \dots, f_{n_jj}$ that correspond to $r_j$. Then, in the Jordan canonical form of $T(\mu)$ there are $n_j$ Jordan blocks whose
diagonal entries are all equal to $r_j$ and the sizes of these Jordan blocks are, say,
$s_{1j} \times s_{1j}, \dots s_{n_jj} \times s_{n_jj}$ (where the $i$-th block
corresponds to the eigenvector $f_{ij}$). Using (\ref{A9aa}) and straightforward induction we get that, for $k = 0, 1, \dots$,
the generalized Floquet solutions (in other words the generalized eigenvectors) $g^l_{ij}$ corresponding to $f_{ij}$ satisfy
\begin{equation}
g^l_{ij}(kb; \mu)
= \sum_{\nu=0}^{l-1} \binom{k}{\nu} r_j^{k-\nu} g^{l-\nu}_{ij}(0; \mu),
\qquad l = 2, \dots, s_{ij},
\label{B6}
\end{equation}
with $g^1_{ij} = f_{ij}$ and the convention that $\binom{k}{\nu} = 0$ for
$k < \nu$.

First, let us assume that $n_j = 1$, for all $j = 1, \dots, m$, namely that, for
each $r_j$ there is only one (up to linear independence) Floquet solution
$f_{1j}(x; \mu)$ such that $f_{1j}(x+b; \mu) = r_j f_{1j}(x; \mu)$ (in other
words we do not have coexistence for any $r_j$). Then, the functions
\begin{equation*}
f_{11}(x), g^2_{11}(x), \dots, g^{s_{11}}_{11}(x); \
f_{12}(x), \dots, g^{s_{12}}_{12}(x); \ \dots; \
f_{1m}(x), \dots, g^{s_{1m}}_{1m}(x)
\end{equation*}
(where $\mu$ is suppressed for typographical convenience) form a basis of
${\cal V}_{\mu}$ (i.e. they are $n$ linearly independent solutions of
$\mathcal{L} u = \mu u$). Thus, $H(\mu) = 0$ is equivalent to
\begin{equation}
\left\vert
\begin{array}{ccccccc}
f_{11}(0) & \cdots & g^{s_{11}}_{11}(0)
& \cdots & f_{1m}(0) & \cdots & g^{s_{1m}}_{1m}(0) \\
f_{11}(b) & \cdots & g^{s_{11}}_{11}(b)
& \cdots & f_{1m}(b) & \cdots & g^{s_{1m}}_{1m}(b) \\
\vdots & \ddots & \vdots & \ddots & \vdots & \ddots & \vdots  \\
f_{11}((n-1)b) & \cdots & g^{s_{11}}_{11}((n-1)b)
& \cdots & f_{1m}((n-1)b) & \cdots & g^{s_{1m}}_{1m}((n-1)b)
\end{array}
\right\vert = 0.
\label{B7}
\end{equation}
In view of (\ref{B6}), the determinant appearing in (\ref{B7}) is equal to
\begin{align}
& f_{11}(0)^{s_{11}} \cdots f_{1m}(0)^{s_{1m}}
\nonumber
\\
\nonumber
\\
& \times \left\vert
\begin{array}{cccccccc}
1 & 0 & \cdots & 0 & \cdots & 1 & \cdots & 0 \\
r_1 & 1 & \cdots & 0 & \cdots & r_m & \cdots & 0 \\
r_1^2 & 2r_1 & \cdots & 0 & \cdots & r_m^2 & \cdots & 0 \\
\vdots & \vdots & \ddots & \vdots & \ddots & \vdots & \ddots & \vdots  \\
r_1^{n-2} & (n-2) r_1^{n-3} & \cdots
& \binom{n-2}{s_{11} - 1} r_1^{n-1 - s_{11}}
& \cdots & r_m^{n-2} & \cdots & \binom{n-2}{s_{1m} - 1} r_m^{n-1 - s_{1m}} \\
\\
r_1^{n-1} & (n-1) r_1^{n-2} & \cdots
& \binom{n-1}{s_{11} - 1} r_1^{n - s_{11}}
& \cdots & r_m^{n-1} & \cdots & \binom{n-1}{s_{1m} - 1} r_m^{n - s_{1m}}
\end{array}
\right\vert .
\nonumber
\\
\label{B8}
\end{align}
By Proposition  A1 of the Appendix the above determinant does not vanish.
Therefore, (\ref{B7}) is equivalent to
\begin{equation*}
f_{11}(0)^{s_{11}} \cdots f_{1m}(0)^{s_{1m}} = 0.
\end{equation*}
Thus, at least one of $f_{11}(0), \dots, f_{1m}(0)$ is $0$, and hence, without
loss of generality we can conclude that there is a $p$ with $1 \leq p \leq m$
such that
\begin{equation}
f_{11}(0) = \cdots = f_{1p}(0) = 0
\qquad \text{while} \qquad
f_{1,p+1}(0) \cdots f_{1m}(0) \ne 0.
\label{B9}
\end{equation}
Since $f_{1j}(x; \mu)$, $j = 1, \dots, m$ are pure Floquet solutions, (\ref{B9})
implies
\begin{equation}
f_{11}(kb) = \cdots = f_{1p}(kb) = 0
\qquad \text{for all }\; k \in \mathbb{Z},
\label{B10}
\end{equation}
hence $f_{1j}(x; \mu)$, $j = 1, \dots, p$ are eigenfunctions of (\ref{B1}) with corresponding eigenvalue $\mu$.

Now, suppose there is another eigenfunction, say $\psi(x)$ corresponding to $\mu$, such that $\psi(x)$ and $f_{1j}(x; \mu)$, $j = 1, \dots, p$ are
linearly independent. Then, without loss of generality $\psi(x)$ must have the form
\begin{align}
\psi(x) = \, & c^1_{1,p+1} f_{1,p+1}(x) + \cdots + c^1_{1m} f_{1m}(x) + c^2_{11} g^2_{11}(x) + \cdots + c^{s_{11}}_{11} g^{s_{11}}_{11}(x)
+ \cdots
\nonumber \\
+ \, & c^2_{1m} g^2_{1m}(x) + \cdots + c^{s_{1m}}_{1m} g^{s_{1m}}_{1m}(x),
\label{B11}
\end{align}
where, of course, at least one of the constants $c^l_{1,j}$ ($l$ is a superscript, as usual) must be nonzero.
Since $\psi(0) = \psi(b) = \cdots = \psi((n-1)b) = 0$, we must have that the columns of the $n \times (n-p)$ matrix
\begin{equation}
\left[
\begin{array}{ccccccc}
g^2_{11}(0) & \cdots & g^{s_{11}}_{11}(0)
& \cdots & f_{1m}(0) & \cdots & g^{s_{1m}}_{1m}(0) \\
g^2_{11}(b) & \cdots & g^{s_{11}}_{11}(b)
& \cdots & f_{1m}(b) & \cdots & g^{s_{1m}}_{1m}(b) \\
\vdots & \ddots & \vdots & \ddots & \vdots & \ddots & \vdots  \\
g^2_{11}((n-1)b) & \cdots & g^{s_{11}}_{11}((n-1)b)
& \cdots & f_{1m}((n-1)b) & \cdots & g^{s_{1m}}_{1m}((n-1)b)
\end{array}
\right]
\label{B12}
\end{equation}
are linearly dependent. Notice that the matrix of (\ref{B12}) arises from the $n \times n$ matrix whose determinant appears in (\ref{B7}),
after erasing the $p$ columns
\begin{equation*}
\left[
\begin{array}{c}
f_{1j}(0) \\
f_{1j}(b) \\
\vdots   \\
f_{1j}((p-1)b)
\end{array}
\right],
\qquad
j = 1, \dots, p.
\end{equation*}
Using (\ref{B6}), the fact (see (\ref{B9})) that $f_{1,p+1}(0) \cdots f_{1m}(0) \ne 0$, and Proposition A1 we can conclude that
\begin{equation*}
g^2_{1j}(0) = 0
\qquad
\text{for some }\; j = 1, \dots, p.
\end{equation*}
Then, (\ref{A9aa}) and (\ref{B9}) imply that
\begin{equation*}
g^2_{1j}(kb) = 0
\qquad
\text{for all }\; k \in \mathbb{Z},
\end{equation*}
thus, these $g^2_{1j}(x)$'s are in ${\cal D}_{\mu}$. If there is, yet, another eigenfunction $\psi(x)$ corresponding to $\mu$,
then the same argument can be applied to show that $g^3_{1j}(x; \mu)$ is also in ${\cal D}_{\mu}$ for some $j$ such that
$g^2_{1j}(x; \mu)$ is in ${\cal D}_{\mu}$, and so on. This completes the proof of the theorem for the case where
$n_j = 1$ for all $j = 1, \dots, m$.

Next, suppose that $n_j > 1$ for some $j = 1, \dots, m$, say $n_1 > 1$. If $f_{11}(0) = f_{21}(0) = \cdots = f_{n_1 1}(0) = 0$,
then all $f_{11},  f_{21}, \dots, f_{n_1 1}$ are in ${\cal D}_{\mu}$. Otherwise, by taking suitable linear combinations of
$f_{11},  f_{21}, \dots, f_{n_1 1}$ we can obtain the (linearly independent) pure Floquet solutions,
say $\hat{f}_{11},  \hat{f}_{21}, \dots, \hat{f}_{n_1, 1}$,
associated to $r_1$ with $\hat{f}_{11}(0) = \hat{f}_{21}(0) = \cdots = \hat{f}_{n_1-1,1}(0) = 0$ and $\hat{f}_{n_1 1}(0) \ne 0$.
Thus $\hat{f}_{i1}$, $i = 1, \dots, n_1-1$, are in ${\cal D}_{\mu}$.

To continue, we look at the generalized Floquet solutions of rank $2$ (see Definition 2) associated to $r_1$. There are at most $n_1$ linearly
independent such solutions. Using linear combinations of these solutions we can make generalized Floquet solutions $g_{i1}^2$ of rank $2$
satisfying (\ref{B1}). Each such solution must satisfy (use (\ref{B6}) with $l=2$)
\begin{equation}
0 = g^2_{i1}(kb) = r_1^k g^2_{i1}(0) + k r_1^{k-1} g^1_{i1}(0),
\qquad
k = 0, 1, \dots, n-1,
\label{B6a}
\end{equation}
where $g^1_{i1}(x)$ is a pure Floquet solution associated to $r_1$. It is straightforward to see that in order to satisfy (\ref{B6a})
we need to have $g^1_{i1}(0) = 0$, from which it follows that
\begin{equation}
g^2_{i1}(kb) = g^1_{i1}(kb) = 0
\qquad \text{for all }\;
k \in \mathbb{Z}.
\label{B6b}
\end{equation}
It follows that, for each generalized Floquet solutions $g_{i1}^2$ of rank $2$ associated to $r_1(\mu)$ which is in ${\cal D}_{\mu}$ there is
a Floquet solutions $g_{i1}^1$ of rank $1$ (i.e. a pure Floquet solution) which is also in ${\cal D}_{\mu}$.

Then, we look for generalized Floquet solutions of rank $3, 4, \dots$, associated to $r_1$, which are in ${\cal D}_{\mu}$. In the same way we
can show that if a generalized Floquet solution $g^l$ of rank $l$ associated to $r_1$ is in ${\cal D}_{\mu}$,
then the correspondind (to $g^l$) Floquet solution of rank $(l-1)$ must also be in ${\cal D}_{\mu}$. We have, thus, determined all pure
and generalized Floquet solutions associated to $r_1$, that are in ${\cal D}_{\mu}$.

Now let $\psi(x)$ be another eigenfunction of (\ref{B1}), which is linearly independent of the above. If $\psi(x)$ is a linear combination only
of pure and generalized Floquet solutions associated to $r_1$, then by using (\ref{B6}) we can show that $\psi(x)$ must be a pure or
generalized Floquet solutions associated to $r_1$ already found, contradicting the assumption that $\psi(x)$ is linearly independent of
those already found. Hence, $\psi(x)$ must be a linear combination of pure and generalized Floquet solutions associated to at least one
additional $r_j$ ($j \ne 1$), not only to $r_1$. We can, then argue as in the previous case of the proof to conclude that there is a pure
Floquet solution associated to some $r_j$ ($j \ne 1$) which is an eigenfunction of (\ref{B1}). We then continue in the same way until we exhaust
all the eigenfunctions of (\ref{B1}) as pure Floquet solution associated to some $r_j$, satisfying the statement of the theorem.
\hfill $\blacksquare$

\smallskip

\textbf{Corollary 4.} If $\psi(x)$ is an eigenfunction of (\ref{B1}), then
\begin{equation*}
\psi(kb) = 0
\qquad
\text{for all }\; k \in \mathbb{Z}.
\end{equation*}

\smallskip

\textbf{Remark 7.} The generalized Floquet solutions in the statement of Theorem 3 cannot be avoided in general. For example, consider the
fourth-order operator $\mathcal{L} u = u'''' + 2a^2 u''$, where $a > 0$ is a constant. Then, if we take $b = \pi / a$, we have that
$\mu = -a^4$ is an eigenvalue of the associated problem (\ref{B1}) and the corresponding eigenspace ${\cal D}_{-a^4}$ is spanned by the solutions $\psi_1(x) = \sin(ax)$ and $\psi_2(x) = x \sin(ax)$ (obviously, the first is a pure Floquet solution, while the second is a generalized Floquet
solution). Only in the special case $m_a(\mu) = 1$ we can be sure that the corresponding eigenspace $\mathcal{D}_{\mu}$ is
(one-dimensional and) spanned by a pure Floquet solution.
\hfill $\diamondsuit$

\smallskip

In the unperturbed case, using (\ref{A6aa}) in formula (\ref{B3}) yields (recalling, also, (\ref{A9b}) for the second equation)
\begin{equation}
\tilde{H}(\lambda) =
\prod_{0 \leq j < k \leq n-1}
\frac{\exp\left(\rho^k \lambda^{1/n} b\right) - \exp\left(\rho^j \lambda^{1/n} b\right)}{\left(\rho^k -\rho^j \right) \lambda^{1/n}}
= \prod_{0 \leq j < k \leq n-1}
\frac{\tilde{r}_{k+1}(\lambda) - \tilde{r}_{j+1}(\lambda)}{\left(\rho^k -\rho^j \right) \lambda^{1/n}},
\label{B4}
\end{equation}
where, as usual, $\lambda^{1/n}$ is the principal $n$-th root of $\lambda$ and $\rho = e^{2 \pi i /n}$.
In particular,
\begin{equation}
\tilde{H}(0) = b^{n (n-1)/2} \ne 0.
\label{B5}
\end{equation}
It follows that $\mu$ is a zero of $\tilde{H}(\lambda)$ if and only if
\begin{equation}
\tilde{r}_j(\mu) = \tilde{r}_k(\mu)
\qquad \text{for some }\;
j,k \in \{1, \dots, n\}, \   j < k
\label{B5aa}
\end{equation}
(thus $\mu = \tilde{\lambda}_{m,d}$ for some $(m, d)$---see (\ref{A9d})).
Then, the associated eigenfunction(s) of (\ref{B1}), for the unperturbed operator $\tilde{\mathcal{L}}$ are
\begin{equation}
\psi_{j,k}(x; \mu) := \exp\left( \rho^{k-1} \mu^{1/n} x \right) - \exp\left( \rho^{j-1} \mu^{1/n} x \right)
\label{B5bb}
\end{equation}
for  $j$ and $k$ satisfying (\ref{B5aa}), where $\mu^{1/n}$ is the principal $n$-th root of $\mu$. Notice that $\psi_{j,k}(x; \mu)$ of (\ref{B5bb})
is a pure Floquet solution corresponding to the multiplier $\tilde{r}_j(\mu) = \tilde{r}_k(\mu)$. Clearly, there is no other (linearly independent)
Floquet solution \textit{corresponding to} $\tilde{r}_j(\mu)$, which satisfies the multipoint boundary conditions of (\ref{B1}). Furthermore,
\begin{equation*}
\frac{d}{d \lambda} \left[ \tilde{r}_k(\lambda) - \tilde{r}_j(\lambda) \right]_{\lambda = \mu}
= \frac{b (\rho^{k-1} - \rho^{j-1})}{n \mu^{(n-1)/n}} \, \tilde{r}_j(\mu) \ne 0.
%\label{B5cc}
\end{equation*}
Therefore, the algebraic and geometric multiplicities of $\mu$ are the same, i.e. in the unperturbed case we have
\begin{equation}
\tilde{m}_a(\mu) = \tilde{m}_g(\mu)
\qquad \text{for all }\;
\mu \in \mathbb{C}.
\label{B5dd}
\end{equation}

We continue by noticing that by (\ref{A9b}) and (\ref{A11}) we have that the discriminant of $\tilde{P}(r; \lambda)$ is
\begin{equation}
D_{\tilde{P}}(\lambda)
= \prod_{0 \leq j < k \leq n-1} \left[\exp\left(\rho^k \lambda^{1 / n} b\right) - \exp\left(\rho^j \lambda^{1 / n} b\right)\right]^2.
\label{B5a}
\end{equation}
Comparison of (\ref{B4}) and (\ref{B5a}) yields
\begin{equation}
D_{\tilde{P}}(\lambda) = (-1)^{(n-1)(n-2)/2} n^n \lambda^{n-1} \tilde{H}(\lambda)^2.
\label{B5b}
\end{equation}
Next, by using the factorization
\begin{equation*}
\frac{e^{a_2 z} - e^{a_1 z} }{(a_2 - a_1) z} = e^{\frac{a_2 + a_1}{2} z}
\prod_{m=1}^{\infty} \left[1 + \frac{(a_2 - a_1) i z}{2 \pi m} \right] \left[1 - \frac{(a_2 - a_1) i z}{2 \pi m} \right]
\end{equation*}
in (\ref{B4}) we obtain
\begin{equation*}
\tilde{H}(\lambda) = b^{n (n-1)/2}
\prod_{m=1}^{\infty} \prod_{d=1}^{n-1} \prod_{j=0}^{n-1} \left[1 - \frac{(\rho^d - 1) \rho^j i \lambda^{1 / n} b}{2 \pi m} \right]
\end{equation*}
or (since $\rho = e^{2 \pi i / n}$ is a primitive $n$-th root of $1$)
\begin{equation}
\tilde{H}(\lambda) = b^{n (n-1)/2}
\prod_{m=1}^{\infty} \prod_{d=1}^{n-1} \left[1 - \frac{i^n (\rho^d - 1)^n b^n \lambda }{2^n \pi^n m^n} \right]
= b^{n (n-1)/2}
\prod_{m=1}^{\infty} \prod_{d=1}^{n-1} \left[1 - \frac{(-1)^d \sin\left(\frac{\pi d}{n}\right)^n b^n \lambda }{\pi^n m^n} \right],
\label{B4a}
\end{equation}
or by recalling (\ref{A9d})
\begin{equation}
\tilde{H}(\lambda) = b^{n (n-1)/2}
\prod_{m=1}^{\infty} \prod_{d=1}^{n-1} \left(1 - \frac{\lambda }{\tilde{\lambda}_{m,d}} \right).
\label{B4d}
\end{equation}
In particular, $\tilde{H}(\lambda)$ has order $1/n$ and all its zeros are real.

In the case where $n$ is even, say $n = 2\nu$, formula (\ref{B4a}) becomes
\begin{equation}
\tilde{H}(\lambda) = b^{(2\nu-1)\nu}
\prod_{m=1}^{\infty} \left\{\left(1 - \frac{(-1)^{\nu} b^{2\nu} \lambda }{\pi^{2\nu} m^{2\nu}} \right)
 \prod_{d=1}^{\nu-1} \left[1 - \frac{(-1)^d \sin\left(\frac{\pi d}{2\nu}\right)^{2\nu} b^{2\nu} \lambda }{\pi^{2\nu} m^{2\nu}} \right]^2
\right\},
\label{B4b}
\end{equation}
while, if $n$ is odd, formula (\ref{B4a}) becomes
\begin{equation}
\tilde{H}(\lambda) = b^{n (n-1)/2}
\prod_{m=1}^{\infty} \prod_{d=1}^{(n-1)/2} \left[1 - \frac{\sin\left(\frac{\pi d}{n}\right)^{2n} b^{2n} \lambda^2 }{\pi^{2n} m^{2n}} \right].
\label{B4c}
\end{equation}

\textbf{Example 2.} (i) If $n = 2$, formula (\ref{B4a}) gives
\begin{equation*}
\tilde{H}(\lambda) = b \prod_{m=1}^{\infty} \left(1 + \frac{b^2 \lambda }{\pi^2 m^2} \right).
\end{equation*}
Here, the set of zeros of $\tilde{H}(\lambda)$ is the Dirichlet spectrum of $d^2 / dx^2$ on $(0, b)$.

(ii) If $n = 3$, formulas (\ref{B4a}) and (\ref{B4c}) give
\begin{equation*}
\tilde{H}(\lambda) = b^3
\prod_{m=1}^{\infty} \left(1 - \frac{3 \sqrt{3} \, b^3 \lambda }{8 \pi^3 m^3} \right)
\left(1 + \frac{3 \sqrt{3} \, b^3 \lambda }{8 \pi^3 m^3} \right)
= b^3 \prod_{m=1}^{\infty} \left(1 - \frac{27 b^6 \lambda^2 }{64 \, \pi^6 m^6} \right).
\end{equation*}
(iii) If $n = 4$, formula (\ref{B4b}) gives
\begin{equation*}
\tilde{H}(\lambda) = b^6
\prod_{m=1}^{\infty} \left(1 - \frac{b^4 \lambda }{\pi^4 m^4} \right)
\left(1 + \frac{b^4 \lambda }{4 \pi^4 m^4} \right)^2.
\end{equation*}
For instance, if $\mu = -4 b^{-4} \pi^4 m^4$ with $m \in \mathbb{N}$, then $\mu$ is a double zero of the $\tilde{H}(\lambda)$
above and its associated eigenfunctions are
\begin{equation*}
e^{\pi m x /b} \sin\left(\frac{\pi m x}{b}\right)
\qquad \text{and} \qquad
e^{-\pi m x /b} \sin\left(\frac{\pi m x}{b}\right).
\end{equation*}
(iv) If $n = 6$, formula (\ref{B4b}) becomes
\begin{equation*}
\tilde{H}(\lambda) = b^{15}
\prod_{m=1}^{\infty}
\left(1 + \frac{b^6 \lambda }{\pi^6 m^6} \right) \left(1 + \frac{b^6 \lambda }{64 \, \pi^6 m^6} \right)^2
\left(1 - \frac{27 b^6 \lambda }{64 \, \pi^6 m^6} \right)^2.
\end{equation*}
Here, let us mention that, for every $k \in \mathbb{N}$ the number $\mu = -64 \, \pi^6 k^6 / b^6$ is a triple zero of $\tilde{H}(\lambda)$
and the associated eigenfunctions are
\begin{equation*}
e^{\sqrt{3}\,\pi k x /b} \sin\left(\frac{\pi k x}{b}\right),
\quad
e^{-\sqrt{3}\,\pi k x /b} \sin\left(\frac{\pi k x}{b}\right),
\quad \text{and} \quad
\sin\left(\frac{2 \pi k x}{b}\right).
\end{equation*}

\textbf{Remark 8.} Let $H_{\xi}(\lambda)$ be the function $H(\lambda)$ of (\ref{B2}) corresponding to the operator $\mathcal{L}_{\xi}$
(see (\ref{TA6})). As we have seen in Subsection 2.5 the characteristic polynomial $P(r; \lambda)$ associated to $\mathcal{L}_{\xi}$ is
independent of $\xi$. Hence, the set $\mathcal{Z}$ (see Definition 1) is independent of $\xi$ too. If $\lambda \not\in \mathcal{Z}$,
then the equation $\mathcal{L}_{\xi} u = \lambda u$ has $n$ linearly independent (pure) Floquet solutions corresponding to $n$ different
Floquet multipliers. It follows that, as $\xi$ moves, the zeros of $H_{\xi}(\lambda)$ which are not in $\mathcal{Z}$ also move, simply
because they are zeros of Floquet solutions. As
for the zeros of $H_{\xi}(\lambda)$ which are in $\mathcal{Z}$, some of them (counting multiplicities) may stay fixed as $\xi$ varies
(for example, if $\mathcal{L}$ is a Hill operator whose spectrum has a ``closed" gap, then $H_{\xi}(\lambda)$ has a simple zero at this
gap for every $\xi$).
\hfill $\diamondsuit$

\subsection{Asymptotics of $H(\lambda)$ as $\lambda \to \infty$}

By using (\ref{A6m}) and (\ref{A6aa}) in (\ref{B3}) one establishes the following theorem (formula (\ref{B4})
explains why we need to exclude the $\lambda$'s such that $\lambda^{1/n} \in E_{\delta}$).

\smallskip

\textbf{Theorem 4.} Let $H(\lambda; t)$, $t \in \mathbb{T}$, be the function defined in (\ref{B2}), corresponding to the
operator $\mathcal{L}(t)$ (see (\ref{TA1})). Then,
\begin{equation}
H(\lambda; t) = \tilde{H}(\lambda)
\left[1 + O\left(\frac{1}{\lambda^{1/n}}\right)\right],
\qquad
\lambda \to \infty,\ \lambda^{1/n} \not\in E_{\delta},
\label{B17}
\end{equation}
uniformly in $t$.

\smallskip

\textbf{Corollary 5.} Given a $\delta > 0$ there is a $K > 0$ such that if a zero $\lambda^{\ast}$ of $H(\lambda; t)$ satisfies
$|\lambda^{\ast}| > K$, then $(\lambda^{\ast})^{1/n} \in E_{\delta}$, where $E_{\delta}$ is the set defined in (\ref{A9g}). The
constants $\delta$ and $K$ are independent of $t$.

\smallskip

The above corollary combined with a theorem of Hurwitz (see, e.g., \cite{H}, Th. 14.3.4) has a remarkable consequence: All zeros of
$H(\lambda; t)$ are limits, as $s \to t$, of zeros of $H(\lambda; s)$, since Corollary 5 and Remark 2 guarantee that no zero of $H(\lambda; s)$
can ``escape" to infinity as $s \to t$.
In other words, as we move $t$ continuously from $t=0$ to $t=1$, the zeros of $H(\lambda; t)$ move continuously in the complex plane
(without escaping to infinity) and hence, each zero of $H(\lambda) = H(\lambda; 1)$ has "evolved" from a zero of
$\tilde{H}(\lambda) = H(\lambda; 0)$ counting multiplicities (and every zero of $\tilde{H}(\lambda)$ evolves to a zero of $H(\lambda)$).

Let us, also, notice that, in particular, Theorem 4, Corollary 5, and the aforementioned ``conservation of the zeros" are valid for the family of
operators $\mathcal{L}_{\xi}$, $\xi \in \mathbb{R}$, introduced in Subsection 2.5.

\smallskip

%\smallskip
%Hence the asymptotics of the fundamental solutions,
%as $\lambda \to \infty$, imples that $H(\lambda)$ is of order $1/n$.??????????????????????????????
%The asymptotic behavior of $H(\lambda)$ as $|\lambda| \to \infty$, is determined by $\tilde{H}(\lambda)$ \cite{N}. One consequence
%is that $H(\lambda)$ has order $1 / n$.

In the self-adjoint Hill (i.e. second-order) case it is easy to see that if $\mu$ is a zero of $H(\lambda)$,
then $m_g(\mu) = m_g(\mu) = 1$. Also, it has been
established in \cite{P1} respectively that $m_a(\mu) = m_g(\mu)$ holds in the case of the fourth-order (periodic) Euler-Bernoulli operator.
We have also seen (recall (\ref{B5dd})) that in the unperturbed case we have equality of the algebraic and the geometric multiplicities.
However, even for the second-order operator $\mathcal{L} = d^2/dx^2 + p_0(x)$ with a nonreal coefficient $p_0(x)$,
it is easy to see that we can have $m_a(\mu) > 1$. For example,
pick any $b$-periodic (nonconstant) entire function $q(z)$ and set $p_0(x) = q(x + \xi)$. Then, it is not hard to see that there are
$\xi \in \mathbb{C}$ for which for the operator $\mathcal{L}$ there is a $\mu$ such that $m_a(\mu) > 1$ (clearly $m_g(\mu) \leq 1$).

One question here is whether $m_a(\mu) = m_g(\mu)$ is always true in the case where $\mathcal{L}$ is a general self-adjoint operator.

%If the coefficients of $\mathcal{L}$ are singular, one can easily construct counterexamples for the above conjecture. For instance, let
%us consider the second-order equation
%\begin{equation}
%\mathcal{L} u = -u'' + \frac{3}{1 - e^{ix}} u = \lambda u.
%\label{Co1}
%\end{equation}
%Here the period is $b = 2\pi$ (the coefficient has simple poles at $x = 2 n \pi$, $n \in \mathbb{Z}$). Observe that the functions
%\begin{equation*}
%f(x) = e^{ix} - e^{2ix}
%\end{equation*}
%and
%\begin{equation*}
%g(x) = 6 e^{ix} \ln(e^{2ix} - 1) + e^{-ix} +3 - 6 e^{ix} - 6 i x e^{ix} -6 e^{2ix} \ln(e^{2ix} - 1) + 6 i x e^{2ix}
%\label{Co2}
%\end{equation*}
%satisfy (\ref{Co1}) for $\lambda = 4$. In fact $f(x)$ is a periodic solution (i.e. a Floquet solution with multiplier $r = 1$), while $g(x)$ is
%???????????????????????????. Also, $f(2 n \pi) = 0$, for all $n \in \mathbb{Z}$. Obviously, there
%is no other (linearly independent) solution that satisfies these boundary conditions. However, since
%\begin{equation}
%\int_0^{2\pi} f(x)^2 dx = 0,
%\label{Co3}
%\end{equation}
%it follows that $\lambda = 4$ is a \textit{double zero} of $H(\lambda)$.

\section{Two conjectures and an open question}

\subsection{A conjecture regarding the structure of the Floquet matrix}

In this subsection the term \textit{lambda matrix} is used for a matrix whose elements are entire functions of $\lambda$
(notice that in the literature the term lambda matrix is sometimes used in a more restrictive way,
namely it indicates that all the elements of the matrix are polynomials in $\lambda$ \cite{M}, \cite{G-K}). The Floquet matrix $T(\lambda)$
is a typical example of a lambda matrix.

Let $N = N(\lambda)$ be an $n \times n$ lambda matrix with characteristic polynomial
\begin{equation}
P_N(r) = P_N(r; \lambda) := \det[r I - N(\lambda)]
\label{D1}
\end{equation}
($I$ is, of course, the $n \times n$ identity matrix). We assume (although it is not necessary) that $P_N(r)$ is irreducible over
$\mathcal{M}(\mathbb{C})$, the field of meromorphic functions on $\mathbb{C}$. In particular, our assumption implies that the discriminant
$D(\lambda)$ of $P_N(r)$ is not identically zero (notice that $D(\lambda)$ is entire in $\lambda$).

The eigenvalues $r_1(\lambda), \dots, r_n(\lambda)$ of $N$ can be viewed as
branches of a multi-valued analytic function $r(\lambda)$. To make $r(\lambda)$ single-valued we consider its
$n$-sheeted Riemann surface $\Sigma$, where the sheets of $\Sigma$ are copies of $\mathbb{C}$.

%The surface $\Sigma$ can be identified with the (generically transcendental, if the elements of $N(\lambda)$ are not polynomials in
%$\lambda$) complex curve over the $\lambda$-plane
%\begin{equation}
%P_N(r; \lambda) = 0.
%\label{D2}
%\end{equation}
If $\lambda$ is not a zero of $D(\lambda)$, then the numbers $r_1(\lambda), \dots, r_n(\lambda)$ are distinct and $N(\lambda)$ is
diagonalizable having $n$ linearly independent eigenvectors, say, $v_1(\lambda), \dots, v_n(\lambda)$, where $v_j(\lambda)$ is associated
to the eigenvalue $r_j(\lambda)$, $j = 1, \dots, n$.

Suppose now that $\ell$ is a ramification point of $\Sigma$ of degree $d > 1$. Then, there are $d$ sheets of $\Sigma$ meeting at $\ell$ corresponding to the branches, say, $r_{j_1}, \dots, r_{j_d}$ of $r$, whose common value at $\ell$ is $r(\ell)$. One might expect that the associated
eigenvectors $v_{j_1}(\lambda), \dots, v_{j_d}(\lambda)$ collapse, as $\lambda$ approaches the projection $\lambda^{\star}$ of $\ell$ on
the complex plane, to one pure eigenvector associated to $r(\ell)$ (i.e. the $d$-dimensional span
$\langle v_{j_1}(\lambda), \dots, v_{j_d}(\lambda) \rangle$
collapses to an eigenline of $N(\lambda^{\star})$), causing the Jordan canonical form of $N(\lambda^{\star})$
to have a $d \times d$ Jordan block. However, this is not always the case. It may happen that
$\langle v_{j_1}(\lambda), \dots, v_{j_d}(\lambda) \rangle$
converges, as $\lambda \to \lambda^{\star}$, to an eigenspace of $N(\lambda^{\star})$ corresponding to $r(\ell)$ of dimension greater than one.
In this case we say that $N(\lambda)$ has a \textit{pathology of the first kind over} $\lambda^{\star}$.

\smallskip

\textbf{Example 3.} Let
\begin{equation*}
N(\lambda) = \left[
\begin{array}{cc}
\lambda^3 + 1 & \lambda^2  \\
2\lambda + \lambda^4 & \lambda^3 + 1  \\
\end{array}
\right].
\end{equation*}
Here,
\begin{equation*}
r(\lambda) = \lambda^3 + 1 + \lambda \sqrt{\lambda (\lambda^3 + 2)},
\end{equation*}
so that $0$ is ramified ($d=2$). However, $N(0)$ is the identity matrix, therefore it has two linearly independent eigenvectors.
Thus, $N$ has a pathology of the first kind over $0$.

\smallskip

\textbf{Conjecture 1.} The Floquet matrix $T(\lambda)$ of (\ref{A3}) does not have pathologies of the first kind over any $\lambda \in \mathbb{C}$.

\smallskip

The conjecture says that if $\ell$ is a ramification point of $\Sigma$ of degree $d > 1$ and $\lambda^{\star}$ is the projection of $\ell$ on the
complex plane, then the Jordan canonical form of the Floquet matrix $T(\lambda^{\star})$ has a $d \times d$ Jordan block associated to the
ramification point $\ell$. In particular, $\lambda^{\star} \in \mathbb{C}$ is not a projection of a ramification
point of $\Gamma$ if $T(\lambda^{\star})$ is diagonalizable (this can happen, of course, even
if $r_j(\lambda^{\star}) = r_k(\lambda^{\star})$ for two different branches $r_j$ and $r_k$ of $r(\lambda)$---e.g., in the
"coexistence" case in the Hill operator theory). If Conjecture 1 is true, the set $\mathcal{R}_{\Gamma}$ is a subset of $\mathcal{Z}_J$ 
(recall Definition 1).

It is well known that Conjecture 1 is true in the Hill case \cite{M-W}. Also, for the fourth-order (periodic) Euler-Bernoulli case Conjecture 1
has been established in \cite{P1}.

Next, suppose that $r(\ell_1) = r(\ell_2) = r^{\star}$, where $\ell_1 \ne \ell_2$ are two points of $\Sigma$ with the same projection,
say $\lambda^{\star}$ on $\mathbb{C}$ (each of $\ell_1$, $\ell_2$ may or may not be ramified). Suppose the sheet of
$r_{j_1}(\lambda)$ passes through $\ell_1$, the sheet of $r_{j_2}(\lambda)$ passes through $\ell_2$ (of course, $j_1 \ne j_2$), and
$r_{j_1}(\lambda^{\star}) = r_{j_2}(\lambda^{\star}) = r^{\star}$. Again, one might expect that the two-dimensional span
$\langle v_{j_1}(\lambda), v_{j_2}(\lambda) \rangle$ of the eigenvectors
$v_{j_1}(\lambda)$ and $v_{j_2}(\lambda)$, associated to $r_{j_1}(\lambda)$ and $r_{j_2}(\lambda)$ respectively, approaches a two-dimensional
eigenspace of $N(\lambda^{\star})$ as $\lambda \to \lambda^{\star}$. However, it may happen that, as $\lambda \to \lambda^{\star}$,
$\langle v_{j_1}(\lambda), v_{j_2}(\lambda) \rangle$ collapses to an one-dimensional eigenline (contributing to the arising of a Jordan anomaly
in $N(\lambda^{\star})$), in which case we say that $N(\lambda)$ has a \textit{pathology of the second kind over} $\lambda^{\star}$.

\smallskip

\textbf{Example 4.} Let
\begin{equation*}
N(\lambda) = \left[
\begin{array}{cc}
(\lambda - 1)^2 & \lambda^2 - 2  \\
\lambda^2 & (\lambda + 1)^2  \\
\end{array}
\right].
\end{equation*}
Here,
\begin{equation*}
r(\lambda) = \lambda^2 + 1 + \lambda \sqrt{\lambda^2 + 2},
\end{equation*}
hence $0$ is not ramified, i.e. $\Sigma$ has two different points whose projection on $\mathbb{C}$ is $0$ (the sheets through these points, over
$0$, do not meet). However, $N(0)$ is not diagonalizable. Thus, $N$ has a pathology of the second kind over $0$.

\smallskip

Notice that a lambda matrix $N(\lambda)$ can have both kinds of pathologies over the same point $\lambda^{\star} \in \mathbb{C}$.
Furthermore, suppose that
$d_1, \dots, d_m$ are the ramification degrees of the points of $\Sigma$ (the Riemann surface associated to $N$) whose projection is
$\lambda^{\star}$, and $s_1, \dots, s_k$ are the sizes of the Jordan blocks of the Jordan canonical form of $N(\lambda^{\star})$. If $N$ has no pathologies over $\lambda^{\star}$, then the sets $\{d_1, \dots, d_m\}$ and $\{s_1, \dots, s_k\}$ are equal (in particular $m = k$).

It is well known that in the self-adjoint Hill case \cite{M-W} the Floquet matrix cannot have pathologies of the second kind. Also, for the
fourth-order (periodic) Euler-Bernoulli case this property has been established in \cite{P1}. Furthermore, the the Floquet matrix of
 the unperturbed $n$-th order cannot have pathologies of any kind. 
However, even for the second-order operator $\mathcal{L} = d^2/dx^2 + p_0(x)$,
one can come up with nonreal coefficients $p_0(x)$ for which the associated Floquet matrix has pathologies of the second kind.
Hence, one may ask whether there are no such pathologies in the case where $\mathcal{L}$ is a general self-adjoint operator.

\subsection{A conjecture regarding the poles of the normalized Floquet solutions}

It has been suggested \cite{D-K-N} and \cite{D-M-N} that the Riemann surface $\Gamma$ and the poles of the normalized Floquet solutions
(as points on $\Gamma$) determine the operator $\mathcal{L}$. This is the inverse periodic spectral problem. Actually, the reconstruction of
$\mathcal{L}$ can be outlined via Abel's theorem and Jacobi inversion \cite{S2}, \cite{S1}. In order for this scheme to work, one
implicitly assumes \cite{No} the following:

\smallskip

\textbf{Conjecture 2.} The poles of the normalized Floquet solution $\phi(x; \lambda)$ are in 1-1
correspondence with the handles of $\Gamma$. In particular, if $\Gamma$ has finitely many ramification points, and hence
it can be compactified to a surface $\overline{\Gamma}$ of finite genus $g$ (see Remark 4),
then $\phi(x; \lambda)$ has exactly $g$ poles on $\Gamma$ (counting multiplicities).

\smallskip

It is well known that this Conjecture is valid in the Hill case. Actually, the projections of these poles are the Dirichlet
eigenvalues lying in (the closure) of nondegenerate gaps in an 1-1 fashion. Also, for the fourth-order (periodic) Euler-Bernoulli case the
Conjecture has been established in \cite{P1}.

\smallskip

\textbf{Remark 9.} If $n = 2$ (i.e. if we are in the Hill case), then the converse of Proposition 1
(see Subsection 2.2) is true (see, e.g., \cite{B} or \cite{M-W}) and sometimes is referred as Borg's Theorem. However, in general the
(naive) converse of Proposition 1 is false (see, e.g., \cite{P} or \cite{P3}). Here is a more carefully stated "converse" of Proposition 1
(even for the non-self-adjoint case), inspired by Remark 3: ``Suppose that the $n$-sheeted Riemann surface $\Gamma$ associated to $\mathcal{L}$ of
(\ref{A1}) has finitely many ramification points and (see Remark 4) its compactification $\overline{\Gamma}$ has genus $0$. Then, the coefficients of $\mathcal{L}$ are constant." For the case $n=3$ this statement has been proved in \cite{Mc}, while for the periodic Euler-Bernoulli operator,
the statement has been proved in \cite{P2}.

If we assume the truth of Conjecture 2, then we can prove the above statement for any order $n$.
Here is how the proof goes: The compact $n$-sheeted surface $\overline{\Gamma}$ is the graph of an algebraic function, and since its genus
is $0$ it must be the graph of a polynomial
\begin{equation}
\lambda = Q(w),
\qquad \text{where }\;
\deg Q = n.
\label{D2}
\end{equation}
In this case Conjecture 2 implies that the normalized Floquet solution $\phi(x; \lambda)$ has no poles on $\Gamma$. Therefore, it
must be an entire function in the variable $w$ (see (\ref{NF7}) and (\ref{NF10})). If $\phi(\xi; \ell) = 0$ for some $\xi \in \mathbb{R}$
and $\ell \in \Gamma$, then $\ell$ must be a pole of the normalized Floquet solution associated to the shifted operator $\mathcal{L}_{\xi}$
(see Subsection 2.5). However, as we have seen the Riemann surface associated to $\mathcal{L}_{\xi}$ is $\Gamma$. Hence, $\phi(x; \lambda)$
cannot vanish. Also, using (\ref{NF7}), (\ref{NF10}), and the fact that the fundamental solutions are entire in $\lambda$ of order $1/n$,
we conclude that $\phi(x; \lambda)$ is entire in $w$ of order one. Therefore, the Hadamard factorization theorem implies that $\phi(x; \lambda)$
must have the form
\begin{equation}
\phi(x; \lambda) = e^{A(x) + B(x) w},
\qquad \text{for some $C^n$ functions }\;
A(x) \text{ and } B(x).
\label{D3}
\end{equation}
Since $\phi(0; \lambda) = 1$, it follows that $B(0) = 0$ (while $A(0)$ must be an integral multiple of $2 \pi i$ and, without loss of generality
we can take $A(0) = 0$). Finally, by using again (\ref{NF7}), (\ref{NF10}), the asymptotic formula (\ref{A6m}), and the fundamental solutions
of the unperturbed case (see (\ref{A6a}) an (\ref{A6aa})), we can conclude that $B(x) = c x$ and $A(x) \equiv 0$, which, in turn, implies
that all the coefficients of $\mathcal{L}$ are constant.
\hfill $\diamondsuit$

\subsection{An open question}

\textbf{Open Question.} Suppose that there exist a polynomial $R(w)$ of degree $\nu$ such that the equation $P(r; \lambda) = 0$ of the multiplier
curve of the operator $\mathcal{L}$ (recall \eqref{A7} and \eqref{A14}) can be written equivalently as
\begin{equation*}
\left[r + r^{-1} - \Delta\left(w_1(\lambda)\right)\right] \cdots \left[r + r^{-1} - \Delta\left(w_{\nu}(\lambda)\right)\right] = 0,
\end{equation*}
where $\Delta(\cdot)$ is an entire function and $w_1(\lambda), \dots, w_{\nu}(\lambda)$ are the solutions of the (polynomial) equation
\begin{equation*}
R(w) - \lambda = 0,
\end{equation*}
i.e. the branches of the $\nu$-valued function $w = R^{-1}(\lambda)$.

Is it, then, true that
\begin{equation*}
\mathcal{L} = R(\mathcal{H}),
\end{equation*}
where $\mathcal{H}$ is a Hill operator? If this is the case, then the function $\Delta$
must be the Hill discriminant of $\mathcal{H}$.

\smallskip

The above question has an affirmative answer in the case of a (fourth-order) periodic Euler-Bernoulli operator \cite{P2}.

\section{Appendix}

\textbf{A1. A Lemma}

\medskip

\textbf{Lemma A1.} Let $z$ be a complex variable and $\sigma > 0$ a fixed constant.

(i) If $|\Re\{z\}| \geq \sigma$, then
\begin{equation}
\left| e^z - 1 \right| \geq 1 - e^{-\sigma}.
\label{AAP0}
\end{equation}
(ii) If $|\Re\{z\}| < \sigma$, then
\begin{equation}
\left| e^z - 1 \right| \geq |\sin R \,| \, \mathbf{1}_{\{R \, \in \mathcal{J}\}} + \mathbf{1}_{\{R \, \in \mathcal{J}^c\}}
+  O\left( \frac{1}{R} \right),
\qquad
R \to \infty,
\label{AAP1}
\end{equation}
where $R := |z|$,
\begin{equation}
\mathcal{J} := \bigcup_{m \in \mathbb{Z}} \left(2 \pi m - (\pi / 2), \, 2 \pi m + (\pi / 2)\right),
\label{AAP2}
\end{equation}
while $\mathbf{1}_S$ denotes the indicator function of the set $S$.

\smallskip

\textit{Proof}. Part (i) is immediate. To prove Part (ii), let us write $z$ in its polar form, namely $z = R e^{i\theta}$.
Since $z \to \infty$ and  $-\sigma < \Re\{z\} < \sigma$ we can assume that $\theta$ satisfies
$|(\pi/2) - \theta| \leq \sigma R^{-1} + O(R^{-3})$ or $|(3\pi/2) - \theta| \leq \sigma R^{-1} + O(R^{-3})$.

Suppose $|(\pi/2) - \theta| \leq \sigma R^{-1} + O(R^{-3})$. Then,
\begin{equation}
|\cos\theta| \leq \frac{\sigma}{R} + O\left( \frac{1}{R^3} \right)
\quad \text{and} \quad
\sin\theta = 1 + O\left( \frac{1}{R^2} \right),
\qquad
R \to \infty.
\label{AAP4}
\end{equation}
Formula (\ref{AAP4}) implies that (as $R \to \infty$)
\begin{equation}
e^z = e^{R \cos \theta} e^{iR \sin \theta} = e^{R \cos \theta} e^{iR}\left[1 + O\left( \frac{1}{R} \right)\right]
= e^{R \cos \theta} e^{iR} + O\left( \frac{1}{R} \right).
\label{AAP5}
\end{equation}
Thus,
\begin{equation}
\left| e^z - 1 \right| = \left| e^{R \cos \theta} e^{iR} - 1\right| + O\left( \frac{1}{R} \right).
\label{AAP6}
\end{equation}
Now, let us consider the quantity.
\begin{equation}
g(x) := \left|x e^{iR} - 1 \right| = \sqrt{x^2 - 2 x \cos R + 1},
\qquad
x \in [0, \infty).
\label{AAP7}
\end{equation}
Notice that $g(x)$ attains its (global) minimum on $[0, \infty)$ at $x^{\ast} := \max\{\cos R, 0\}$. Hence
\begin{equation}
g(x) \geq g(x^{\ast}) =  |\sin R \,| \, \mathbf{1}_{\{\cos R \; > \, 0\}} + \mathbf{1}_{\{\cos R \; \leq \, 0\}}
\qquad \text{for all }\; x \in [0, \infty).
\label{AAP8}
\end{equation}
Using (\ref{AAP7}) and (\ref{AAP8}) in (\ref{AAP6}) yields
\begin{equation}
\left| e^z - 1 \right|  \geq |\sin R \,| \, \mathbf{1}_{\{R \, \in \mathcal{J}\}} + \mathbf{1}_{\{R \, \in \mathcal{J}^c\}}
+ O\left( \frac{1}{R} \right),
\qquad
R \to \infty.
\label{AAP9}
\end{equation}
The case $|(3\pi/2) - \theta| \leq \sigma R^{-1} + O(R^{-3})$ is reduced to the case $|(\pi/2) - \theta| \leq \sigma R^{-1} + O(R^{-3})$
by noticing that $|e^z - 1| = |e^{\bar{z}} - 1|$.
\hfill $\blacksquare$

\smallskip

\textbf{Corollary A1.} Given $\beta \in (0, 1)$ there is a $\delta = \delta(\beta) > 0$ such that
\begin{equation*}
\left| e^z - 1 \right| \geq \beta
\qquad \text{for all }\;
z \in \bigcup_{m \in \mathbb{Z}} \{z \in \mathbb{C}:\, |z - 2 \pi i m| \geq \delta\}.
\end{equation*}
Furthermore, $\delta = \delta(\beta)$ can be chosen so that $\delta \to 0^+$ as $\beta \to 0^+$.

\bigskip

\textbf{A2. The Gershgorin Circle Theorem}

\medskip

Here we remind the reader of a basic theorem of linear algebra.

\smallskip

\textbf{Theorem A1} (the Gershgorin Circle Theorem)\textbf{.} Let $A = [a_{jk}]$ be an $n \times n$ matrix. Then, the eigenvalues of $A$
(and hence of $A^{\top}$, since $A^{\top}$ and $A$ have the same characteristic polynomial) are located in the union of the $n$ disks
(the so-called \textit{Gershgorin disks}):
\begin{equation}
\left\{r \in \mathbb{C}:\, |r - a_{jj}| \leq \sum_{k \ne j} |a_{jk}| \right\},
\qquad
j = 1, \dots, n.
\label{TA13}
\end{equation}
Furthermore, if the $j$-th Gershgorin disk is isolated from the other disks, then it contains exactly one eigenvalue of $A$.

\smallskip

The idea of the proof is simple \cite{S-S}: For each eigenvalue $r$ of $A$ we normalize the associated eigenvector $x = (x_1, \dots, x_n)^{\top}$
so that its component with the maximum absolute value, say $x_j$, is set equal to $1$ (if there are more than one components with maximum
absolute value, just choose one of them and set it equal to $1$). Then, $r x = A x$ implies
\begin{equation*}
r = r x_j = \sum_{k=1}^n a_{jk} x_k = a_{jj} + \sum_{k \ne j} a_{jk} x_k,
\end{equation*}
from which (\ref{TA13}) follows immediately, since $|x_k| \leq 1$ for all $k = 1, \dots, n$.

The last statement of the theorem follows by viewing $A$ as a continuous deformation of the matrix $A_0 := \text{diagonal} (a_{11}, \dots, a_{nn})$
and then using the fact that every eigenvalue of a matrix is a continuous function of the elements of the matrix.

\bigskip

\textbf{A3. Generalized (or confluent) Vandermonde determinants}

\medskip

Let $A \geq 1$ and $\alpha \geq 1$ be integers with $A \geq \alpha$. We introduce the $A \times \alpha$ block (i.e. matrix)
\begin{equation}
B(x; A \times \alpha) = (c_{jk})_{1 \leq j \leq A \atop 1 \leq k \leq \alpha} \, ,
\qquad \text{where} \quad
c_{jk} := \binom{j-1}{k-1} x^{j - k}
\label{AP0}
\end{equation}
(recall the convention $\binom{j-1}{k-1} = 0$ for $j < k$). Notice that $B(x; A \times \alpha)$ is a square matrix only if $A = \alpha$, and in
this case its determinant is $1$.

Next, let $\alpha_1, \dots, \alpha_m$ be $m$ strictly positive integers and
\begin{equation}
A = \alpha_1 + \cdots + \alpha_m \, .
\label{AP1}
\end{equation}
Putting the blocks $B(x_1; A \times \alpha_1), \dots, B(x_m; A \times \alpha_m)$ side by side we form the $A \times A$ (square) matrix
\begin{equation}
M(x_1, \dots, x_m; \alpha_1, \dots, \alpha_m) := \left[ B(x_1; A \times \alpha_1) \cdots B(x_m; A \times \alpha_m) \right]
\label{AP0a}
\end{equation}
and consider its determinant
\begin{equation}
F(x_1, \dots, x_m; \alpha_1, \dots, \alpha_m) := \det M(x_1, \dots, x_m; \alpha_1, \dots, \alpha_m),
\label{AP0b}
\end{equation}
namely
\begin{align}
& F(x_1, \dots, x_m; \alpha_1, \dots, \alpha_m)
\nonumber
\\
\nonumber
\\
& = \left\vert
\begin{array}{cccccccc}
1 & 0 & \cdots & 0 & \cdots & 1 & \cdots & 0 \\
x_1 & 1 & \cdots & 0 & \cdots & x_m & \cdots & 0 \\
x_1^2 & 2x_1 & \cdots & 0 & \cdots & x_m^2 & \cdots & 0 \\
\vdots & \vdots & \ddots & \vdots & \ddots & \vdots & \ddots & \vdots  \\
x_1^{A-2} & (A-2) x_1^{A-3} & \cdots
& \binom{A-2}{\alpha_1 - 1} x_1^{A-1 - \alpha_1}
& \cdots & x_m^{A-2} & \cdots & \binom{A-2}{\alpha_m - 1} x_m^{A-1 - \alpha_m} \\
\\
x_1^{A-1} & (A-1) x_1^{A-2} & \cdots
& \binom{A-1}{\alpha_1 - 1} x_1^{A - \alpha_1}
& \cdots & x_m^{A-1} & \cdots & \binom{A-1}{\alpha_m - 1} x_m^{A - \alpha_m}
\end{array}
\right\vert .
\nonumber
\\
\label{AP2}
\end{align}
For instance, if $m = 3$ and $(\alpha_1, \alpha_2, \alpha_3) = (2, 3, 1)$ we get
\begin{equation*}
F(x_1, x_2, x_3; 2, 3, 1)
= \left\vert
\begin{array}{cccccc}
1      & 0        & 1      & 0        & 0          & 1      \\
x_1    & 1        & x_2    & 1        & 0          & x_3    \\
x_1^2  & 2 x_1    & x_2^2  & 2 x_2    & 1          & x_3^2  \\
x_1^3  & 3 x_1^2  & x_2^3  & 3 x_2^2  & 3 x_2      & x_3^3  \\
x_1^4  & 4 x_1^3  & x_2^4  & 4 x_2^3  & 6 x_2^2    & x_3^4  \\
x_1^5  & 5 x_1^4  & x_2^5  & 5 x_2^4  & 10 x_2^3   & x_3^5
\end{array}
\right\vert
= (x_2 - x_1)^6 (x_3 - x_1)^2 (x_3 - x_2)^3.
\end{equation*}
In the case $\alpha_1 = \cdots = \alpha_A = 1$ (hence $m = A$)
$F(x_1, \dots, x_A; 1, \dots, 1)$ becomes the standard Vandermonde determinant
and we have
\begin{equation*}
F(x_1, \dots, x_A; 1, \dots, 1)
= \prod_{1 \leq j < k \leq A} (x_k - x_j).
\end{equation*}
On the other hand, in the extreme case $m = 1$ we have $\alpha_1 = A$ and
\begin{equation*}
F(x_1; A) \equiv 1.
\end{equation*}
\textbf{Proposition A1.} If $m \geq 2$, then
\begin{equation}
F(x_1, \dots, x_m; \alpha_1, \dots, \alpha_m)
= \prod_{1 \leq j < k \leq m} (x_k - x_j)^{\alpha_j \alpha_k}.
\label{AP3a}
\end{equation}

For the proof see \cite{H-J}.


\begin{thebibliography}{4}

\bibitem{A} L.V. Ahlfors, \textit{Complex Analysis}, Third Edition, McGraw-Hill Book Company, New York, 1979.

\bibitem{B-K} A. Badanin and E. Korotyaev, Even order periodic operators on the real line, \textit{IMRN}, \textbf{2012} (No. 5), (2012),
1143--1194. doi:10.1093/imrn/rnr057

\bibitem{B} G. Borg, Eine umkehrung der Sturm-Liouvilleschen eigenwertaufgabe, \textit{Acta Math.}
\textbf{78} (1946), 1--96.

\bibitem{B-K} V.M. Bukhshtaber and I.M. Krichever, Vector Addition Theorems and Baker-Akhiezer Functions,
\textit{Teoreticheskaya i Matematicheskaya Fizika}, \textbf{94} (2) (February 1993), 200--212.
Translation: \textit{Theoretical and Mathematical Physics}, \textbf{94} (2) (1993), 142--149.

%\bibitem{B-C1} J.L. Burchnall and T.W. Chaundy, Commutative ordinary differential operators, \textit{Proc. London Math. Soc. Ser. 2}
%\textbf{21} (1923), 420--440 (DOI: 10.1112/plms/s2-21.1.420).
%
%\bibitem{B-C2} J.L. Burchnall and T.W. Chaundy, Commutative ordinary differential operators, \textit{Proc. Roy. Soc. London A}
%\textbf{118} (1928), 557--583.
%
\bibitem{C-L} E.A. Coddington and N. Levinson, \textit{Theory of Ordinary Differential Equations}, Robert E. Krieger Publishing Company, Malabar,
Florida, 1987.

\bibitem{M0} M.L. Da Silva Menezes, \textit{Infinite genus curves with hyperelliptic ends}, Ph.D. Thesis,
New York University, 1986 (Thesis Advisor: Henry P. McKean, Jr.).

\bibitem{M} M.L. Da Silva Menezes, Infinite genus curves with hyperelliptic ends,
\textit{Comm. Pure Appl. Math.} \textbf{42} (2) (1989), 185--212.

\bibitem{D-K-N} B.A. Dubrovin, I.M. Krichever, and S.P. Novikov, Integrable
Systems, I, Dynamical Systems, IV, 177--332, Encyclopaedia of
Mathematical Sciences, 4, Springer, Berlin, 2001.

\bibitem{D-M-N} B.A. Dubrovin, V. Matveev, and S.P. Novikov, Nonlinear
Equations of Korteweg-deVries Type, Finite Zone Linear Operators, and
Abelian Varieties, \textit{Uspekhi. Mat. Nauk}, \textbf{31} (1976), 55--136; \textit{Russ. Math. Surveys}, \textbf{31} (1976), 59--146.

\bibitem{D-S} N. Dunford and J.T. Schwartz, \textit{Linear Operators. Part II: Spectral Theory; Self Adjoint Operators in Hilbert Space}, Wiley Classics Library Edition, New York, 1988.

\bibitem{E} M.S.P. Eastham, \textit{The spectral theory of periodic differential equations},
Texts Math., Scott., Acad. Press, Edinburgh, 1973.

\bibitem{G-K} I. Gohberg,  M.A. Kaashoek, and L. Lerer, The resultant for regular matrix polynomials and quasi commutativity,
\textit{Indiana Univ. Math. J.}, \textbf{57} (no. 6) (2008), 2793--2813.

\bibitem{H} E. Hille, \textit{Analytic Function Theory, Vol. II}, Second Edition, Chelsea Publishing Company, New York, 1977.

\bibitem{H-J} R. Horn and C. Johnson, \textit{Topics in Matrix Analysis}, Cambridge University Press, 1991.

\bibitem{K} M.V. Keldysh, On the characteristic values and characteristic functions of certain classes of non-self-adjoint
equations, \textit{Doklady Akad. Nauk SSSR} (N.S.) \textbf{77} (1951), 11--14 (Russian); \textit{Math. Rev.}, \textbf{12} 835 (1951).

\bibitem{K1} P. Kuchment, \textit{Floquet Theory for Partial Differential
Equations}, Birkhauser, Verlag, Basel, 1993.

\bibitem{K2} P. Kuchment, private communication (1996).

\bibitem{L} S. Lang, \textit{Algebra}, Sixth Printing, Addison-Wesley Publishing Company, Reading, Massachusetts, 1974.

\bibitem{M-W} W. Magnus and S. Winkler, \textit{Hill's Equation}, Dover Publications, Inc., New York, 1979.

\bibitem{M} A.S. Markus, \textit{Introduction to the spectral theory of polynomial operator pencils}. Translated from the Russian by H. H. McFaden. Translation edited by Ben Silver. With an appendix by M.V. Keldysh. Translations of Mathematical Monographs, \textbf{71}, American Mathematical Society, Providence, RI, 1988.

\bibitem{Mc} H.P. McKean, Boussinesq's equations on the circle, \textit{Comm. Pure Appl. Math.}, \textbf{34} (1981), 599--691.

\bibitem{M-Mo} H.P. McKean and P. Moerbeke, The spectrum of Hill's operator, \textit{Invent. Math.}, \textbf{30} (3) (1975), 217--274.

\bibitem{M-M} H.P. McKean and V.H. Moll, \textit{Elliptic Curves, Function Theory, Geometry, Arithmetic}, Cambridge University
Press, New York, 1999.

\bibitem{M-T} H.P. McKean, and E. Trubowitz, Hill's operator and hyperelliptic function theory in the presence of
infinitely many branch points, \textit{Ann. Pure Appl. Math.} \textbf{29} (1976), 143--226.

\bibitem{N} M.A. Naimark, \textit{Linear Differential Operators}, Parts I \& II, Frederick Ungar Publishing Co., New York, 1967 \& 1968.

\bibitem{No} S.P. Novikov, private communication.

\bibitem{P} V.G. Papanicolaou, The Spectral Theory of the Vibrating Periodic Beam,
\textit{Communications in Mathematical Physics}, \textbf{170} (1995), 359--373.

\bibitem{P1} V.G. Papanicolaou, The Periodic Euler-Bernoulli Equation, \textit{Transactions of the American Mathematical Society}, \textbf{355} (9)
(2003), 3727--3759.

\bibitem{P2} V.G. Papanicolaou, An Inverse Spectral Result for the Periodic Euler-Bernoulli Equation, \textit{Indiana Univ. Math. Journal},
\textbf{53} (1) (2004), 223--242.

\bibitem{P3} V.G. Papanicolaou, The Inverse Periodic Spectral Theory of the Euler-Bernoulli Equation,
\textit{Dynamics of Partial Differential Equations}, \textbf{2} (2) (2005), 127--148.

\bibitem{P4} V.G. Papanicolaou, Generalized Vandermonde Determinants, preprint (2014).

\bibitem{R-S} M. Reed and B. Simon, \textit{Methods of modern mathematical physics. I. Functional analysis}, Second edition. Academic Press, Inc. (Harcourt Brace Jovanovich, Publishers), New York, 1980.

\bibitem{S2} C.L. Siegel, \textit{Topics in Complex Function Theory, Vol. II: Automorphic Functions and Abelian Integrals}, Wiley
Classics Library Edition, Wiley--Interscience, a division of John Wiley and Sons, Inc., New York, 1988.

%\bibitem{S-S-dJ} S. Silvestrov, C. Svensson, and M. de Jeu,
%Algebraic Dependence of Commuting Elements in Algebras
%
%\bibitem{S} B. Simon, The Classical Moment Problem as a Self-Adjoint
%Finite Difference Operator, \textit{Advances in Mathematics} \textbf{137} (1998),
%82--203 (Article No. AI981728).
%
\bibitem{S1} G. Springer, \textit{Introduction to Riemann Surfaces}, 2nd ed., Chelsea, New York, 1981.

\bibitem{S-S} G.W. Stewart and J.-G. Sun, \textit{Matrix perturbation theory}, Academic Press, Boston, 1990.

\bibitem{T} E.C. Titchmarsh, \textit{The Theory of Functions}, Second Edition, Oxford, 1939.

\bibitem{Tk} V.A. Tkachenko, Spectral Analysis of a Nonselfadjoint Hill Operator,
\textit{Dokl. Akad. Nauk SSSR}, \textbf{322} (No. 2) (1992); \textit{Soviet Math. Dokl.}, \textbf{45} (No. 1) (1992), 78--82.

\bibitem{Y-S} V.A. Yakubovich and V.M. Starzhinskii, \textit{Linear differential equations with periodic coefficients. 1, 2}. Translated from
Russian by D. Louvish. Halsted Press (John Wiley \& Sons) New York-Toronto, Ont.,; Israel Program for Scientific Translations, Jerusalem-London, 1975.

\end{thebibliography}
\end{document}